\documentclass[10pt,a4paper]{article}

\usepackage[english]{babel}
\usepackage{amsthm,amsmath,amssymb,amsfonts}
\usepackage[utf8]{inputenc}
\usepackage{graphicx}
\usepackage{hyperref}
\usepackage{authblk}

\newcommand{ \B }{{ \bf b}}
\newcommand{ \BB }{{ \bf B}}
\newcommand{ \C}{ {\bf c}}
\newcommand{ \DD }{{ \bf D}}
\newcommand{ \ctilde }{c_{\sharp}}

\newcommand{\e}{ \mathrm{e} }

\newcommand{\ed}[1]{ a_{#1} }
\newcommand{\F}{ \mathcal{F} }
\newcommand{\I}{ \mathrm{i} }

\newcommand{ \K }{ {\bf k}}
\newcommand{ \Kperp }{ {\bf k}_\perp}
\newcommand{ \N }{\nu_{\sharp}}
\newcommand{ \p }{ {\bf p}}

\newcommand{ \Pdiff }{ \mathcal{P}}
\newcommand{ \q }{ {\bf q}}

\newcommand{ \R }{ {\bf r}}

\newcommand{ \U }{ {\bf u}}

\newcommand{\W}{ \mathcal{W} }

\newcommand{ \x }{ {\bf x}}
\newcommand{ \xperp }{ {\bf x}_\perp}

\newcommand{ \xn }{ x_n}
\newcommand{ \y }{ {\bf y}}
\newcommand{ \yperp }{ {\bf y}_\perp}

\newcommand{ \z }{ {\bf z}}
\newcommand{\Z}{ \mathcal{Z} }

\title{Weak localization in radiative transfer of acoustic waves in a randomly-fluctuating slab}

\author[1,2]{Adel Messaoudi}
\author[1]{R\'egis Cottereau}
\author[2]{Christophe Gomez}
\affil[1]{{Aix Marseille Univ, CNRS, Centrale Marseille, LMA}, Marseille, France}
\affil[2]{{Aix Marseille Univ, CNRS, I2M}, Marseille, France}

\begin{document}
\maketitle
\begin{abstract}
This paper concerns the derivation of radiative transfer equations for acoustic waves propagating in a randomly fluctuating slab (between two parallel planes) in the weak-scattering regime, and the study of boundary effects through an asymptotic analysis of the Wigner transform of the wave solution. These radiative transfer equations allow to model the transport of wave energy density, taking into account the scattering by random heterogeneities. The approach builds on the method of images, where the slab is extended to a full-space, with a periodic map of mechanical properties and a series of sources located along a periodic pattern. Two types of boundary effects, both on the (small) scale of the wavelength, are observed: one at the boundaries of the slab, and one inside the domain. The former impact the entire energy density (coherent as well as incoherent) and is also observed in half-spaces.
The latter, more specific to slabs, corresponds to the constructive interference of waves that have reflected at least twice on the boundaries of the slab and only impacts the coherent part of the energy density.
\end{abstract}

\begin{flushleft}
\end{flushleft}



\section{Introduction}
\label{sec:intro}
Radiative transfer theory was introduced over a century ago to describe the propagation of light in complex media. Today, it is used in many other fields such as in geophysics \cite{Maeda2008,Margerin2019,Margerin1998}, neutronics \cite{Larsen1975,Larsen1976,Larsen1974,Spanier2008}, optics \cite{Klose2010,Klose1999}, for weather forecasting \cite{Villefranque2019}, or for the illumination of animation movies scenes \cite{Bitterli2018,Keller2015}. Radiative transfer equations can be derived from the wave equation in the high-frequency regime \cite{Bal2005,Bal2010,Bal2002,Cottereau2014,butz2015,Ryzhik1996} through a multi-scale asymptotic analysis of the Wigner transform of the wave field. 

In this paper, we more particularly concentrate on (scalar) acoustics in a slab, which is a domain contained in between two parallel planes:
\[
\Omega=\mathbb{R}^2 \times \left(0,H\right).
\]
In this case, the wave equation for the pressure field $p(t,\x)$ is given by:
\begin{equation}
\label{eq:equation_des_ondes}
\partial^2_{tt} p(t,\x)-c^2(\x)\Delta p(t,\x) = 0, \qquad \left(t, \x \right) \in \mathbb{R}_+^*\times \Omega,
\end{equation} 
where $c(\x)$ is the sound speed in the medium, modeled as a random field and assumed to fluctuate at the scale $\ell_c$ (correlation length). Neumann boundary conditions complete this wave equation:
\begin{equation}
\label{eq:neumann_condition}
\partial_{n}p(t,\xperp,\xn=0)=\partial_{n}p(t,\xperp,\xn=H)=0, \qquad \left(t, \xperp \right) \in \mathbb{R}_+^*\times \mathbb{R}^2 ,
\end{equation}
where the spatial variable $\x$ has been split into $\x = (\xperp,\xn)$, and $\xn$ represents the coordinate along the unit outward normal to the interface $\mathbb{R}^2 \times \{0\}$, $\xperp$ its transverse coordinate, and $\partial_n$ stands for the derivative with respect to the variable $\xn$. The case of Dirichlet boundary conditions and mixed boundary conditions 
are considered respectively in Sec.~\ref{sec:dirichlet_boundary_condition} and Sec.~\ref{sec:dirichlet__and_neumann_boundary_conditions}. Finally, the following initial conditions are considered:
\begin{equation}
\label{eq:initial_conditions_wave}
p(t=0,\x)=A\left(  \x-\x_{0}   \right)\quad \mathrm{and} \quad \partial_t p(t=0,\x)=B\left(  \x-\x_{0}   \right), \qquad\x\in\Omega,
\end{equation} 
where the shape of the functions $A(\x)$ and $B(\x)$ defines the wavelength $\lambda$. 

In the so-called weak-scattering regime, where the parameter $\epsilon = \lambda/L$ is small ($L$ being the propagation length), and the velocity field $c(\x)$ fluctuates weakly ($\mu^2\approx\epsilon$) at the scale of the wavelength ($\ell_c\approx\lambda$), the energy density can be shown to verify a radiative transfer equation (RTE):
\begin{equation}
\label{eq:ETR}
\partial_t W +  c_0\widehat{\K}\cdot\nabla_{\x} W =- \Sigma(\K) W  +  \int_{\mathbb{R}^3}\sigma(\K,\q)W(\q) \delta(c_0(|\q|-|\K|))d\q. 
\end{equation}
Here, $W(t,\x,\K)$ represents the energy density at time $t$ and position $\x$ in direction $\K$ (see Sec.~\ref{sec:wigner_transform} for a precise definition of Wigner transform), $\widehat{\K}=\K/|\K|$ is the normalized wave direction, and $c_0$ is the average sound speed in the medium. The differential scattering cross-section $\sigma(\K,\q)$ represents the rate at which energy density along wave vector $\q$ is diffracted into energy density along wave vector $\K$, and the total scattering cross-section is
\begin{equation}
\Sigma(\K) = \int_{\mathbb{R}^3}\sigma(\q,\K)\delta(c_0(|\q|-|\K|))d\q.
\end{equation}
This asymptotic result is usually obtained for problems supported on the entire space ($\Omega=\mathbb{R}^3$), and has been derived starting from various wave equations using different approaches (see references above). The case of media delimited by boundaries has been more marginally studied. The main configuration studied in the literature corresponds to a half-space ($\Omega=\mathbb{R}^2 \times \mathbb{R}^\ast_-$) \cite{Jin2006b,Gerard1993, Ryzhik1997, Burq1997, Messaoudi2022}, for which radiative transfer models are still valid with boundary conditions on $\partial\Omega$ reminiscent of geometrical optics:
\begin{equation}
\label{eq:reflecting_boundary_conditions}
W(t,\left(\xperp,\xn = 0), (\Kperp,k_n)\right) = W(t,\left(\xperp,\xn = 0), (\Kperp,-k_n)\right), 
\end{equation}
where $ \left(t, \xperp \right) \in \mathbb{R}_+^*\times \mathbb{R}^2.$
These results have been extended to rough interfaces~\cite{Bal1999} and some models display coupling between surface and bulk waves~\cite{Margerin2019, Borcea2021}. Moreover, interference effects within one wavelength of the interface have been described~\cite{Messaoudi2022} ($\partial\Omega=\mathbb{R}^2\times \{0\}$), resulting in a doubling of the total energy with respect to the standard radiative transfer model with Neumann boundary conditions (and a cancellation in the case of Dirichlet boundary conditions). 

\begin{figure}[tbhp]
\centering \includegraphics[trim = 1.1cm 1.7cm 0.1cm 0.5cm, clip, scale=0.18]{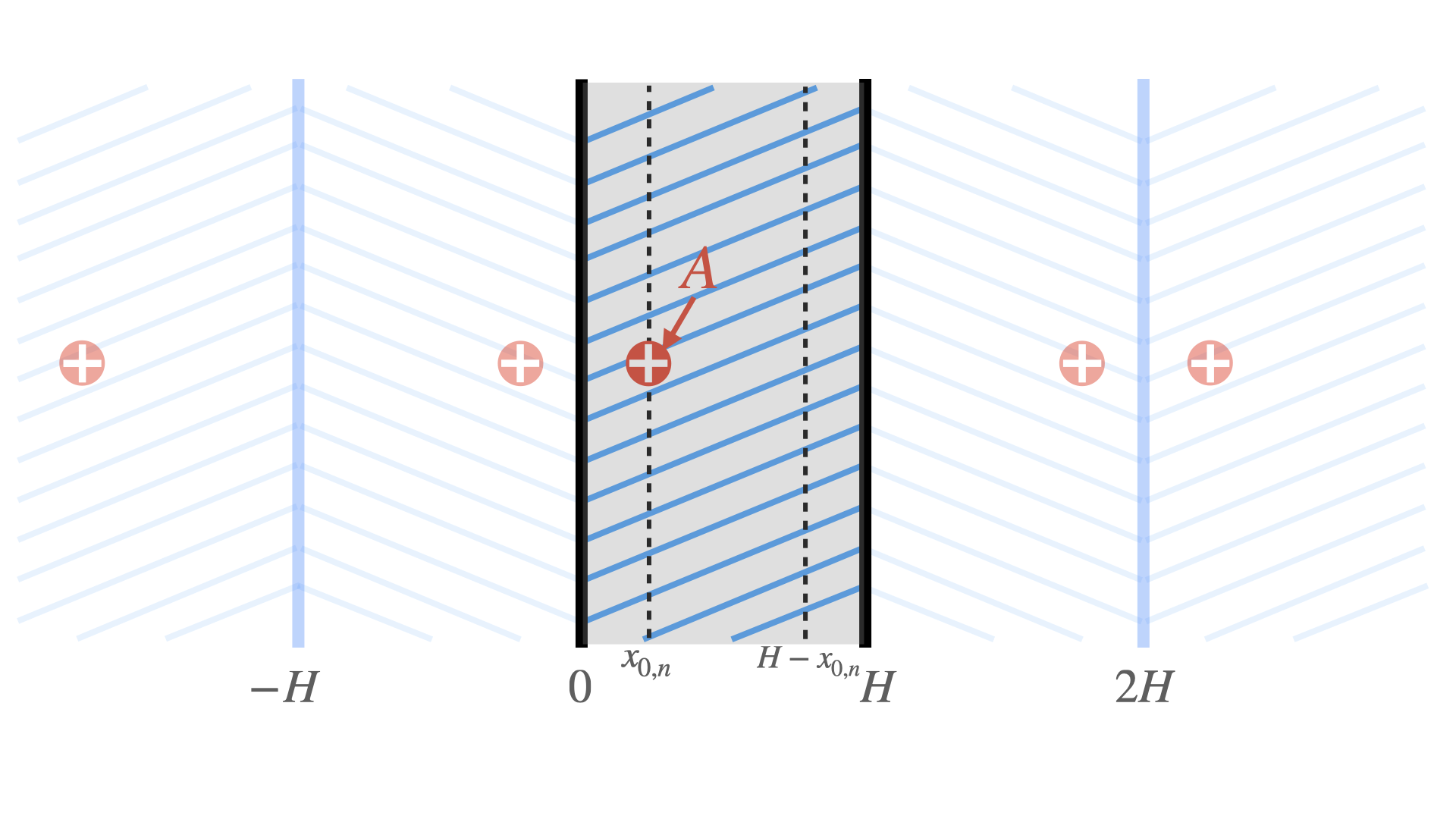}
\caption{\label{fig:superposition_ondes} Sketch of the method of images for the medium delimited by the boundaries $\{\xn=0\}$ and $\{\xn=H\}$. Each wave propagates in a periodic full-space from a boundary condition whose support is centered on one of the red crosses. An amplification of the total energy is observed within one wavelength of the boundaries of the slab (thick black lines), as well as an amplification of the coherent energy inside the slab, along the dashed lines.}
\end{figure}

In this paper, we consider a case that has not been considered before in the context of RTEs, that of acoustic waves propagating in the weak-scattering regime in a randomly-fluctuating slab $\Omega=\mathbb{R}^2 \times \left(0,H\right)$ (see Fig.~\ref{fig:superposition_ondes}), whose thickness $H$ is of the same order as the propagation length $L$ (large compared to the wavelength). As will be shown below, the situation in a slab is very different from the situation in a half-space (and obviously from the full-space). In addition to interference phenomena taking place at the boundaries of the slab in a periodic manner and similar to the ones observed for the half-space, interference will take place between the waves propagating to the left after bouncing on the right-most boundary and the waves propagating to the right after bouncing on the left-most boundary. Similar interference will also take place in a periodic manner between waves bouncing back and forth between the boundaries in different directions. The interference will only take place between the coherent parts of the waves, and will result in amplifications along two different planes of thickness $\epsilon$: one passing through the support of the initial condition and one at the same distance from the other boundary (see the dashed lines in Fig.~\ref{fig:superposition_ondes}). 

As will appear clearly in the mathematical derivations, the propagation in a slab considered in this paper bears similarities with the propagation in a periodic medium that was considered for instance in~\cite{Bal2018}, where RTEs are derived for solutions of Schr\"odinger equations under the weak-scattering regime. However, the period was of the same order as the wavelength, which resulted in major differences (order 1) with respect to RTEs in a full space. On the contrary, the slab we consider in this paper has a thickness of order $H\approx1$, and the interference effects that will be identified take place in areas of thickness $\epsilon$. 

Following the approach discussed in~\cite{Messaoudi2022} for the half-space, our approach is based on the method of images. It consists in representing the solution $p(t,\x)$ to Eq.~\eqref{eq:equation_des_ondes}, supported on $\Omega$, with the superposition of several solutions of the same wave equation \eqref{eq:equation_des_ondes} but extended to the full-space $\mathbb{R}^3$ with periodic sound speed field. In this superposition, one term is equipped with the original initial condition, and the others with initial conditions chosen so as to enforce the boundary conditions (see Fig.~\ref{fig:superposition_ondes} for a sketch, and details in Sec.~\ref{sec:the_method_of_images}). 
This approach allows to obtain the main result of this paper, which is to unveil energy amplification within one wavelength (of order $\epsilon$) of the boundaries (the solid black lines in Fig.~\ref{fig:superposition_ondes}), as well as weak localization phenomena within one wavelength of two planes inside the slab (the dashed lines in Fig.~\ref{fig:superposition_ondes}). The amplifications along the boundaries apply to the entire energy density field while the amplifications inside the domain only apply to the coherent field. Note that, although the terminology \emph{weak localization} (also sometimes called \emph{coherent backscattering}) is mostly used in the literature to refer to constructive interference between wave that have been scattered on the heterogeneities of a half-space (see for instance~\cite{Larose2004}), we use it in this paper in the sense considered in~\cite{Catheline2011,Gallot2011}, of interference effects between coherent waves that have bounced on the boundaries of a domain.

The outline of this paper is as follows. In Sec.~\ref{sec:acoustic_waves_in_a_random_mirror_in_the_high_frequency_regime}, the method of images is introduced, as well as the asymptotic regime under consideration and the Wigner transform. In Sec.~\ref{sec:asymptotic_analysis_with_an_initial_condition_far_from_the_border}, the RTEs are derived using the method of images and interference and weak localization phenomena are highlighted in the case of homogeneous Neumann conditions at both boundaries of the slab. Finally, in Sec.~\ref{sec:other_boundary_conditions}, the influence of considering other sets of boundary conditions is discussed.


\section{Acoustic waves in a randomly fluctuating slab}
\label{sec:acoustic_waves_in_a_random_mirror_in_the_high_frequency_regime}

In the following, we consider the scalar wave equation~\eqref{eq:equation_des_ondes} over the propagation medium $\Omega=\mathbb{R}^2 \times (0,H)$, equipped with the boundary conditions~\eqref{eq:neumann_condition} and the initial conditions~\eqref{eq:initial_conditions_wave}. It is assumed that the initial conditions $A$ and $B$ are smooth functions on $\mathbb{R}^3$, compactly supported within $\Omega$, and that they do not cross the boundaries, so as to be compatible with all the boundary conditions considered in this paper. For simplicity $A$ and $B$ are also assumed to be even w.r.t the $\xn$-variable. Although this assumption is not necessary, and the proposed analysis can be extended to more general situations, it simplifies greatly the forthcoming presentation. We next describe the method of images, the weak-scattering regime, detail how the random inhomogeneities are modeled, and introduce the Wigner transform.


\subsection{The method of images}
\label{sec:the_method_of_images}

The basic principle of the method of images is to replace an acoustic problem posed on a slab, with a given initial condition (or source), by a family of problems on the full-space with a set of initial conditions and a periodic extension of the sound speed, such that the superposition (restricted to the slab) of the solutions associated to the individual problems satisfies the original acoustic problem. 

More precisely, we first consider the following extension $\ctilde(\x)$ of the original velocity field $c(\x)$ to the full-space (see Fig.~\ref{fig:superposition_ondes}):
\begin{multline}
\label{eq:sound_speed_extension}
\ctilde(\x)  = \sum_{j\in\mathbb{Z}} c\left(\xperp,\xn-2jH\right)\mathrm{1}_{[0,H)} (\xn-2jH) \\
+ c\left(\xperp,-\xn+2jH\right)\mathrm{1}_{[-H,0]} (\xn-2jH)  \qquad \x \in \mathbb{R}^3.
\end{multline}
Note that $\ctilde$ is obtained by, first, extending $c$ over $\mathbb{R}^2\times (-H,H)$ by parity (as an even function with respect to the $x_n$-variable), and then over $\mathbb{R}^3$ by periodicity in $x_n$. The resulting function $\ctilde$ is therefore $2H$-periodic and even with respect to the $x_n$-variable. Secondly, we consider the solution (on the full-space $\mathbb{R}^3$), denoted $q$, to the following wave equation associated to the extended sound speed $\ctilde$,  
\begin{equation}
\label{eq:equation_des_onde_full_space}
\partial^2_{tt} q(t,\x)-\ctilde^2\left(\x \right) \Delta q(t,\x) = 0 \qquad  \left(t, \x \right) \in \mathbb{R}_+^*\times \mathbb{R}^3,
\end{equation} 
with the original initial conditions: 
\begin{equation}
\label{eq:IC_p-+}
\begin{cases}
q(t=0,\x) = A\left( \x-\x_{0} \right), \\
\partial_t q (t=0,\x)  = B\left( \x-\x_{0} \right).
\end{cases}
\end{equation}
Finally, we introduce $p^\sharp$,
\begin{equation}
\label{eq:superposition}
p^\sharp(t,\x) :=  \sum_{j\in \mathbb{Z}} q(t,\xperp, x_n - 2jH) + q(t,\xperp, -(x_n - 2jH)), \quad (t,\x) \in  \mathbb{R}_+\times \mathbb{R}^3,
\end{equation}
which satisfies Eq.~\eqref{eq:equation_des_onde_full_space}. The above choices induce the following symmetry relations (evenness and $2H$-periodicity):
\begin{equation}
\label{eq:sym_p-p+}
p^\sharp(t,\xperp,-\xn) = p^\sharp(t,\xperp,\xn) \quad \mathrm{and} \quad  p^\sharp(t,\xperp,\xn + 2jH) = p^\sharp(t,\xperp,\xn),           
\end{equation}
where $(t,\x) \in  \mathbb{R}_+\times \mathbb{R}^3$ and $ j \in\mathbb{Z}$. In particular, this ensures that the derivative $\partial_n p^\sharp(t,\x)$ vanishes at $\xn=0$ and $\xn=H$. As a result, the restriction of $p^\sharp(t,\x)$ to $\Omega$ corresponds to the solution of the original problem:
\begin{equation}
\label{def:p_field}
p^\sharp(t,\x) = p(t,\x), \quad \left(t, \x \right) \in \mathbb{R}_+\times \Omega.
\end{equation}
The above construction corresponds to Neumann boundary conditions applied to the original problem, but other boundary conditions  (homogeneous Dirichlet and mixed boundary conditions) will be considered in detail in Sec.~\ref{sec:other_boundary_conditions}.

To summarize, the solution $p(t,\x)$ of Eq.~\eqref{eq:equation_des_ondes}, defined on the slab $\Omega = \mathbb{R}^2\times (0,H)$, has been extended, through $\ctilde$ and $p^\sharp$, to a periodically-structured full-space. The extension $p^\sharp(t,\x)$ is the superposition of several solutions of the same wave equation~\eqref{eq:superposition} with different initial conditions whose supports are arranged as illustrated in Fig.~\ref{fig:superposition_ondes}. As discussed further down (see Sec.~\ref{sec:wigner_transform}), the description of the overall energy propagation on the slab at leading order only requires to describe the energy propagation of each solution in Eq.~\eqref{eq:superposition} over the full-space. Interesting interference effects will be observed on a small scale (order $\epsilon$) by the interactions of pairs of waves, along the mediator plane of the segment linking the supports of their respective initial conditions. As a result, for all pairs of initial conditions, and because our interest is only on the restriction of $p^\sharp(t,\x)$ to $[0,H]$, these interactions can only take place either along the boundaries of the slab $\{ \xn=0\}$ and $\{\xn=H\}$, or along the planes $\{\xn=x_{0,n}\}$ and $\{\xn=H-x_{0,n}\}$.

\subsection{Weak-scattering regime.}
\label{sec:high_frequency_regime}

The weak-scattering regime corresponds to the scaling $t \to t/\epsilon,$ $\x \to \x/\epsilon$ for a small parameter $\epsilon = \ell_c/L\approx\lambda/L\ll 1$  and $\mu^2 \approx \epsilon$. In other words, the typical wavelength $\lambda$ is of the same order of magnitude as the correlation length $\ell_c$, both are small compared to the typical propagation length $L$ and slab thickness $H$ ($H\approx L$), and the medium velocity fluctuates weakly at the scale of the wavelength. We thus set
\begin{equation*}
p_\epsilon(t,\x) := p\left(\frac{t}{\epsilon},\frac{\x}{\epsilon} \right), \quad \left(t, \x \right) \in \mathbb{R}_+^*\times \Omega,
\end{equation*}
and in a similar fashion the rescaled $q_\epsilon$ associated to $q$. The orders of derivation in time and space being the same for the wave equation, the equilibrium equations~\eqref{eq:equation_des_ondes} and~\eqref{eq:equation_des_onde_full_space}, for $p_\epsilon$ and $q_\epsilon$ are unchanged. The initial conditions are rescaled as
\begin{equation}
\label{eq:cond_init_eps}
p_\epsilon(t=0,\x)=\frac{1}{\epsilon^{3/2}}A\left( \frac{ \x-\x_{0} }{\epsilon}  \right)\quad\text{and}\quad \epsilon \, \partial_t p_\epsilon(t=0,\x)=\frac{1}{\epsilon^{3/2}}B\left(  \frac{ \x-\x_{0} }{\epsilon}  \right),
\end{equation}
where the amplitudes are chosen so as to provide solutions of order 1 for the RTE at the limit $\epsilon\to 0$, and similarly for the initial conditions~\eqref{eq:IC_p-+}. We also consider the rescaled $p^\sharp_\epsilon$, associated to $p^\sharp$, defined by
\begin{equation}
\label{eq:superposition_H.F}
p^\sharp_\epsilon(t,\x) :=  \sum_{j\in \mathbb{Z}} q_\epsilon(t,\xperp, x_n - 2jH) + q_\epsilon(t,\xperp, -x_n + 2jH), \quad (t,\x) \in  \mathbb{R}_+\times \mathbb{R}^3.
\end{equation}
In this paper, the velocity field is assumed to fluctuate around a constant background value $c_0$ as follows (a slowly fluctuating background can also be considered, see for instance~\cite{Ryzhik1996}):
\begin{equation}
\label{eq:sound_speed}
c^2\left( \frac{\x}{\epsilon} \right):= c_0^2 - \sqrt{\epsilon} \, \nu \left( \frac{\x}{\epsilon}\right) \qquad \x\in \Omega,
\end{equation}
where $\nu(\x)$ accounts for the random fluctuations. The amplitude $\sqrt{\epsilon}$ is the proper scaling that provides, in the asymptotic $\epsilon\to 0$, a nontrivial radiative transfer model describing the energy propagation. Under our high-frequency regime, we assume that the extended velocity field $\ctilde$ scales as
\begin{multline}
\label{eq:sound_speed_extension_H.F}
\ctilde\left(\x,\frac{\x}{\epsilon}\right)  = \sum_{j\in\mathbb{Z}}  c\left(\frac{\xperp}{\epsilon},\frac{\xn-2jH}{\epsilon}\right)\mathrm{1}_{[0,H)} (\xn-2jH) \\
+ c\left(\frac{\xperp}{\epsilon},\frac{-\xn+2jH}{\epsilon}\right)\mathrm{1}_{[-H,0]} (\xn-2jH)  \qquad \x \in \mathbb{R}^3.
\end{multline} 
In this latter expression, the slow and fast $\x$-components are separated, the former being used to delimit the periodicity cells. With this choice, the slab width $H$ is of order the typical propagation length $L$ (of order one in the scaling introduced above), which is large compared to the wavelength.


Even though it is possible to work directly with the second-order form of the wave equation to derive the radiative transfer model, we rather use the formalism described in~\cite{Bal2005} based on a first order hyperbolic system. To derive this system from the wave equation, we introduce the vector field 
\begin{equation}
\label{eq:equation_des_ondes_ordre_1}
\U_\epsilon^\sharp(t,\x):=\begin{pmatrix} p^\sharp_\epsilon(t,\x)   \\ \epsilon\,  \ctilde^{-2}\left(\x,\frac{\x}{\epsilon}\right)   \partial_t p^\sharp_\epsilon(t,\x) \end{pmatrix}, \quad (t,\x)\in  \mathbb{R}_+\times\mathbb{R}^3.
\end{equation}
Through the relation~\eqref{eq:superposition_H.F}, $\U_\epsilon^\sharp$ is composed of a sum of vectors,
\begin{equation}
\label{eq:superposition_vecteurs}
\U^\sharp_\epsilon(t,\x)= \sum_{j\in\mathbb{Z}} \U_\epsilon(t,\xperp,x_n-2jH) + \U_\epsilon(t,\x, -x_n +2jH),
\end{equation}
where 
\begin{equation}
\label{eq:equation_des_ondes_ordre_1bis}
\U_\epsilon(t,\x):=\begin{pmatrix} q_\epsilon(t,\x)   \\ \epsilon\,  \ctilde^{-2}\left(\x,\frac{\x}{\epsilon}\right)   \partial_t q_\epsilon(t,\x) \end{pmatrix}, \qquad (t,\x)\in  \mathbb{R}_+\times\mathbb{R}^3,
\end{equation}
satisfies the $2\times2$ system of equations
\begin{equation}
\label{eq:equation_hyperbolique_ordre_1}
\epsilon\partial_t \U_\epsilon + \mathcal{A}_\epsilon\U_\epsilon = 0, \qquad \mathrm{ with }\qquad \mathcal{A}_\epsilon:=-\begin{pmatrix}0& \ctilde^2 \left(\x, \frac{\x}{\epsilon} \right)\\ \epsilon^2\Delta & 0  \end{pmatrix}, 
\end{equation}
equipped with initial conditions inherited from Eq.~\eqref{eq:IC_p-+}.


\subsection{Structure of the inhomogeneities}
\label{sec:structure_of_the_inhomogeneities}

In Eq.~\eqref{eq:sound_speed}, the random fluctuations of the sound-speed profile are modeled through the random field $\nu(\x)$ restricted to $\Omega$. In this paper, $\nu$ is a statistically homogeneous mean-zero random field, defined on $\mathbb{R}^3$, with (normalized) power spectrum density $\widehat{R}(\p)$ given by
\begin{equation*}
(2\pi)^3c_0^4\widehat{R}(\p)\delta(\p+\q):=\big\langle\widehat{\nu}(\p)\widehat{\nu}(\q)\big\rangle.
\end{equation*}
Here, $\langle\cdot\rangle$ denotes an ensemble average and $\widehat{\nu}(\K)$ refers to the Fourier transform of $\nu(\x)$, with the convention
\begin{equation*}
\widehat{\nu}(\K) := \int_{\mathbb{R}^3} \e^{ - \I \K \cdot \x} \nu(\x) d\x, \qquad \mathrm{and} \qquad \nu(\x) := \frac{1}{(2\pi)^3}\int_{\mathbb{R}^3} \e^{  \I \K \cdot \x}\widehat{\nu}(\K) d\K.
\end{equation*}
For the sake of simplicity in the computations, we assume that $R$, the correlation function of $\nu$, is of the form $R(\x)=r(|\x|)$. This assumption implies on its Fourier transform $\widehat{R}$ the following property:
\begin{equation}
\label{eq:R_hat}
\widehat{R}(\p)=\widehat{R}(\p_\perp,-\p_n)=\widehat{R}(-\p_\perp,\p_n)=\widehat{R}(-\p).
\end{equation}

The random field corresponding to the extension~\eqref{eq:sound_speed_extension_H.F} is given by
\begin{multline}
\label{def:N}
\N\Big(\x,\frac{\x}{\epsilon}\Big):=\sum_{j\in\mathbb{Z}} \nu \left(\frac{\xperp}{\epsilon},\frac{\xn-2jH}{\epsilon}\right)\mathrm{1}_{[0,H)} (\xn-2jH) \\ + \nu \left(\frac{\xperp}{\epsilon},\frac{-\xn+2jH}{\epsilon}\right)\mathrm{1}_{[-H,0]} (\xn-2jH),
\end{multline}
which is not stationary anymore w.r.t. the $x_n$-variable. Using the following representation
\begin{multline}
\N\left(\x,\frac{\x}{\epsilon}\right)=\N(\x,\y)_{|\y = \frac{\x}{\epsilon}}=\sum_{j\in\mathbb{Z}} \left(  \nu\left(\y_\perp,y_n-\frac{2jH}{\epsilon}\right)\mathbf{1}_{[0,H)}(\xn-2jH) \right. \\
 \left.+  \nu \left(\y_\perp,-y_n+\frac{2jH}{\epsilon}\right)\mathbf{1}_{[-H,0]}(\xn-2jH)\right)\, {}_{\big|\y = \frac{\x}{\epsilon}},
\end{multline}
the power spectrum density of $\N(\x,\y)$ with respect to the fast component $\y$ reads
\begin{equation}
\label{eq:power_spectra}
\Big\langle \widehat{\N}(\x,\p)  \widehat{\N}(\x,\q)\Big \rangle= (2\pi)^3c_0^4\delta(\p+\q)\widehat{R}(\p),
\end{equation}
for any fixed slow component $\x$, and where $\widehat{\N}(\x,\p)$ stands for the Fourier transform with respect to the $\y$-variable. Note that the power spectrum associated to $\N(\x,\y)$ does not depend on the slow component $\x$, so that the extension procedure will play no role, at leading order in $\epsilon$, in the limiting RTE described in Sec.~\ref{sec:asymptotic_analysis_with_an_initial_condition_far_from_the_border}.


\subsection{Wigner transform}
\label{sec:wigner_transform}

Under the weak-scattering regime, the derivation of RTE from the wave equation relies on a multi-scale asymptotic analysis of the Wigner transform of the wave field. The Wigner transform of two vector fields $\mathbf{v}$ and $\mathbf{w}$ is defined as
\begin{equation*}
W[\mathbf{v},\mathbf{w}](t,\x,\K):=\int_{\mathbb{R}^3}\e^{\I\K\cdot \y} \mathbf{v}\left(t,\x- \frac{\epsilon\y}{2}\right)\otimes \mathbf{w}\left(t,\x+ \frac{\epsilon \y}{2}\right)\frac{d\y}{(2\pi) ^3},
\end{equation*}
where $\otimes$ stands for the tensor product. We may think of the Wigner transform as the inverse Fourier transform of the two point correlation function of $\mathbf{v}$ and $\mathbf{w}$. The Wigner transform is also often called energy density since it relates to the energy carried by a wave field $\U$ through the relation
\[
\mathbf{v}(t,\x)\otimes \mathbf{v}(t,\x) = \int_{\mathbb{R}^3} W[\mathbf{v},\mathbf{v}](t,\x,\K) d\K.
\]

\begin{figure}[tbhp]
\centering \includegraphics[trim = 1.1cm 1.0cm 0.1cm 0.5cm, clip, scale=0.18]{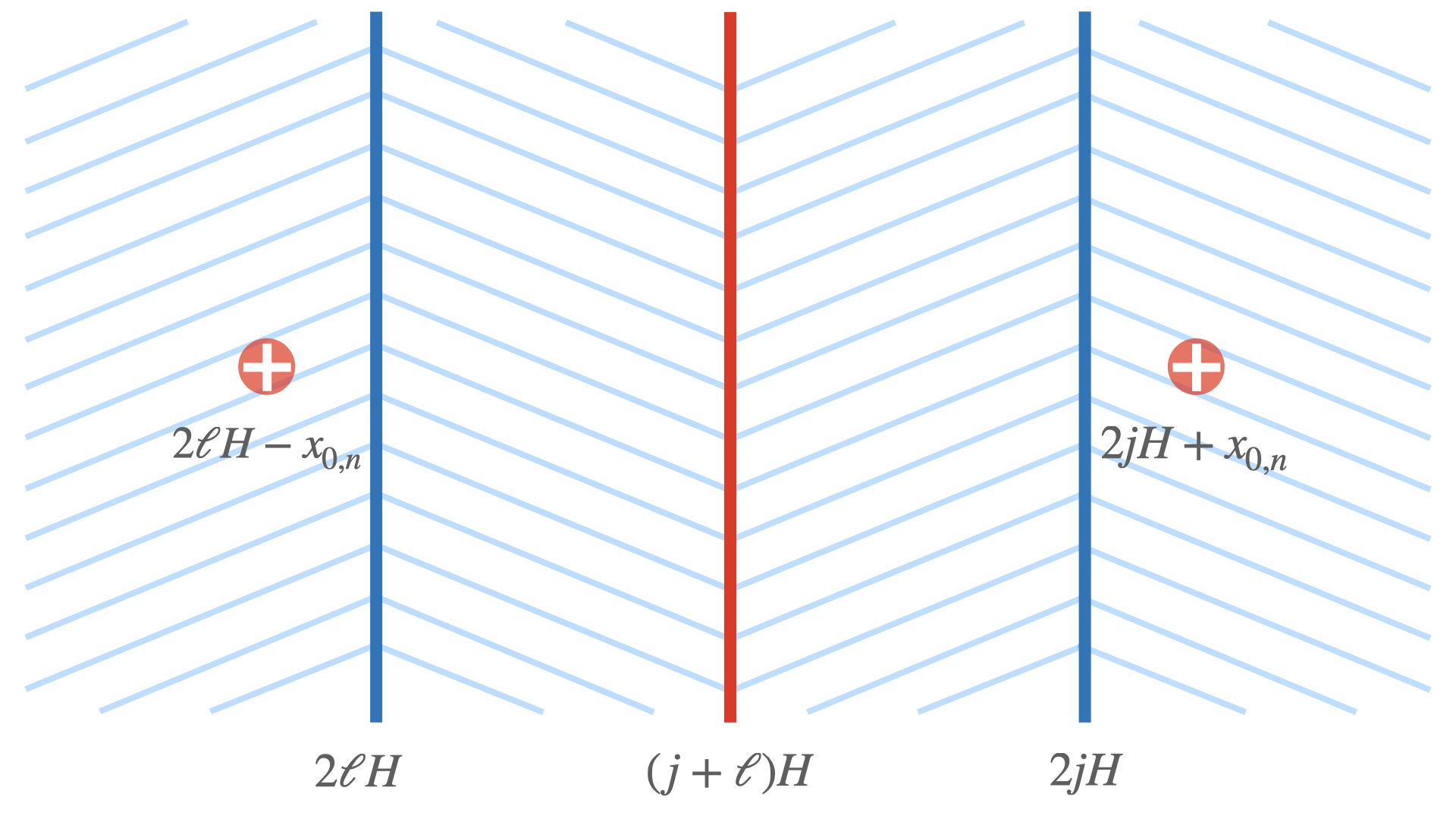}
\caption{\label{fig:cross_wigner1} Each cross-Wigner transform $V_\epsilon^{j\ell}$ corresponds to pairs of waves propagating from initial conditions whose supports are centered around the two red crosses. An amplification is observed along the middle plane at $\xn = (j+\ell)H$ (in red).}
\end{figure}

\begin{figure}[tbhp]
\centering \includegraphics[trim = 1.1cm 1.0cm 0.1cm 0.5cm, clip, scale=0.18]{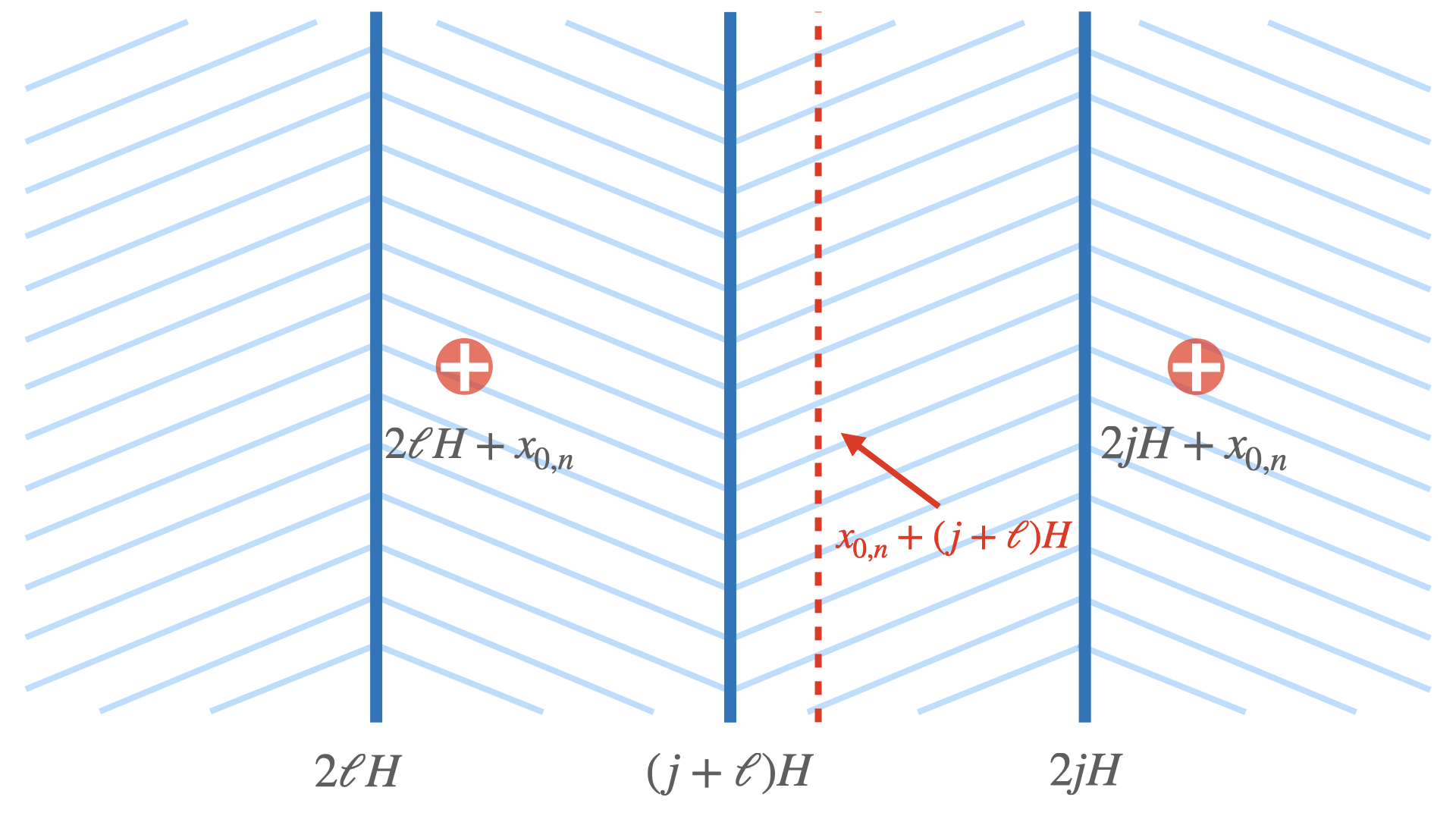}
\centering \includegraphics[trim = 1.1cm 1.0cm 0.1cm 0.5cm, clip, scale=0.18]{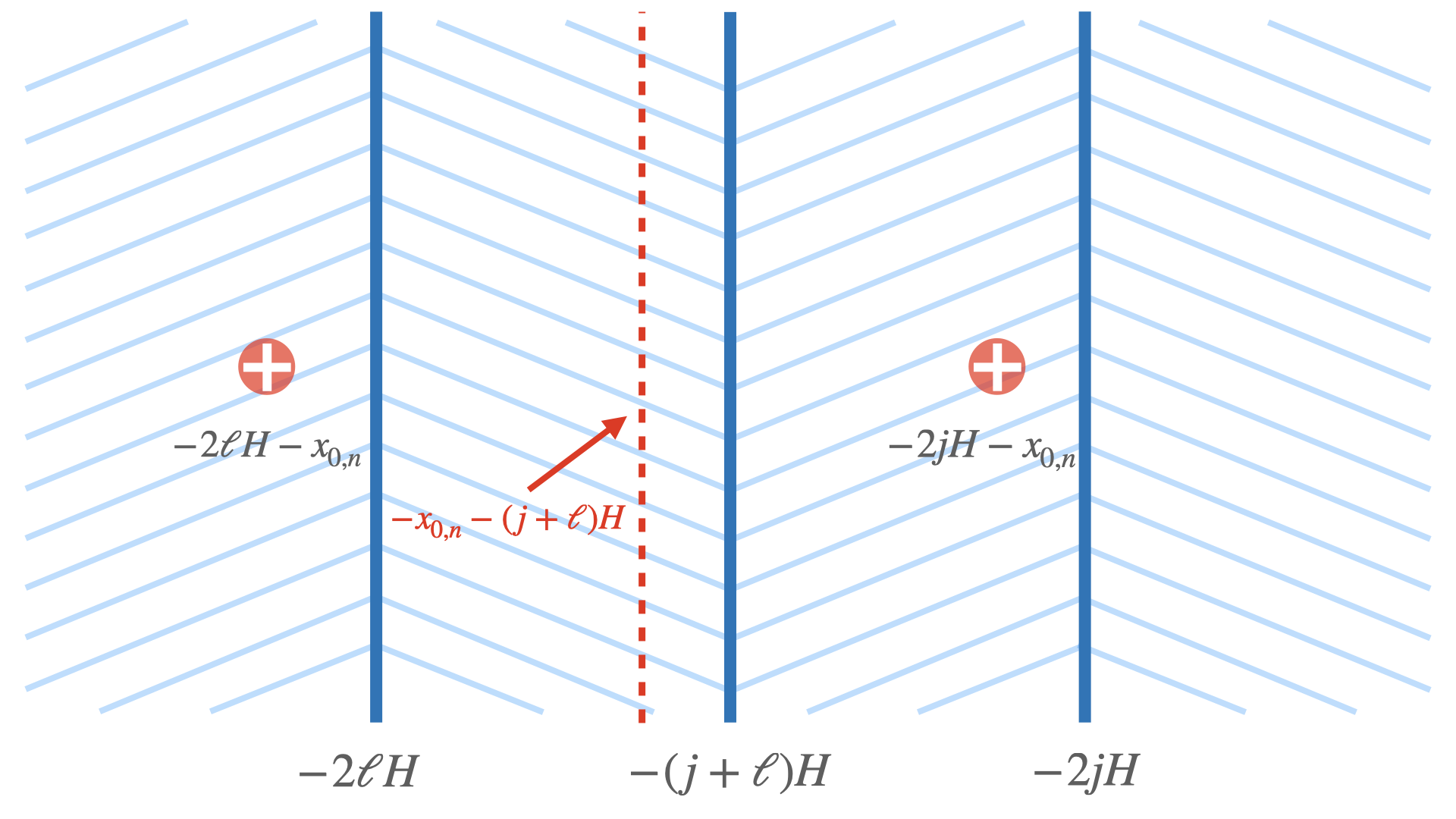}
\caption{\label{fig:cross_wigner3} Each cross-Wigner transform $W_\epsilon^{j\ell}$ corresponds to pairs of waves propagating from initial conditions whose supports are centered around the two red crosses, in two possible configurations. An amplification of the coherent energy is observed along the middle plane at $\xn = \pm x_{0,n} \pm (j+\ell)H$ (dashed red line), depending on the configuration.  } 
\end{figure}

The Wigner transform of $\U_{\epsilon}^\sharp$, defined by Eq.~\eqref{eq:equation_des_ondes_ordre_1} and associated to the extended wave problem, reads    
\begin{equation}
\label{eq:wigner_u}
W_\epsilon^\sharp(t,\x,\K) := W[\U_{\epsilon}^\sharp,\U_{\epsilon}^\sharp](t,\x,\K), \quad \left(t,\x,\K\right)\in\mathbb{R}_+\times\mathbb{R}^3\times\mathbb{R}^3. 
\end{equation}
This Wigner transform~\eqref{eq:wigner_u} can be decomposed, according to Eq.~\eqref{eq:superposition_vecteurs}, as
\begin{multline}\label{eq:wigner_total}
W_\epsilon^\sharp(t,\x,\K)  =  \sum_{j\in \mathbb{Z}} W_\epsilon\left(t,(\xperp,x_n-2jH), \K\right) +  W_\epsilon\left(t,(\xperp,-x_n+2jH),(\Kperp,-k_n) \right)      \\
 +  \sum_{j,\ell \in \mathbb{Z} } V_\epsilon^{j\ell}(t,\x,\K) + V_\epsilon^{j\ell}\left(t,(\xperp,-x_n), (\Kperp,-k_n)\right)\\
 + \sum_{\substack{ j,l \in \mathbb{Z}   \\ \ell \ne j}} W_\epsilon^{j\ell}(t,\x,\K) + W_\epsilon^{j\ell}\left(t,(\xperp,-x_n), (\Kperp,-k_n)\right)
\end{multline}
where three different types of Wigner transform appear: $W_\epsilon$, $V_\epsilon^{j\ell}$ and $W_\epsilon^{j\ell}$, for all $j,\ell\in\mathbb{Z}$. The first type is simply the Wigner transform associated to $\U_{\epsilon}$, and will denoted hereafter as self-Wigner transform:
\begin{equation}
\label{eq:self_wigner}
W_\epsilon(t,\x,\K) :=  W[  \U_\epsilon ,\U_\epsilon] (t,\x,\K)
\end{equation}
The second type $V_\epsilon^{j\ell}(t,\x,\K)$ corresponds to the Wigner transform of a pair of waves whose initial conditions are centered, respectively, around $(\x_{0,\perp},2jH+x_{0,n})$ and $(\x_{0,\perp},2\ell H-x_{0,n})$:
\begin{multline}
V_\epsilon^{j\ell}(t,\x,\K) := \int_{\mathbb{R}^3}\e^{\I\K\cdot \y}  \U_\epsilon\left(t,\left(\xperp- \frac{\epsilon\yperp}{2} ,  \xn - \frac{\epsilon y_n}{2} -2jH \right) \right) \nonumber\\ 
 \otimes \U_\epsilon \left(t,\left(\xperp+ \frac{\epsilon\yperp}{2} ,  -\xn - \frac{\epsilon y_n}{2} +2\ell H \right) \right)\frac{d\y}{(2\pi) ^3}.
\end{multline}
As illustrated in Fig.~\ref{fig:cross_wigner1}, the mediator plane between the supports of the corresponding two initial conditions is located at $\xn=(j+\ell)H$.
Finally, for all $j,\ell \in \mathbb{Z}$, the third type $W_\epsilon^{j\ell}(t,\x,\K)$ corresponds to the Wigner transform of a pair of waves whose initial conditions are centered, respectively, around $(\x_{0,\perp},2jH+x_{0,n})$ and $(\x_{0,\perp},2\ell H+x_{0,n})$:
\begin{multline*}
W_\epsilon^{j\ell}(t,\x,\K) := \int_{\mathbb{R}^3}\e^{\I\K\cdot \y}  \U_\epsilon\left(t,\left(\xperp- \frac{\epsilon\yperp}{2} ,  \xn - \frac{\epsilon y_n}{2} -2jH \right) \right) \\ 
 \otimes \U_\epsilon \left(t,\left(\xperp+ \frac{\epsilon\yperp}{2} ,  \xn + \frac{\epsilon y_n}{2} - 2\ell H \right) \right)\frac{d\y}{(2\pi) ^3}
\end{multline*}
As illustrated in Fig.~\ref{fig:cross_wigner3}, the mediator plane between the supports of the corresponding two initial conditions is located at $\xn= x_{0,n} + (j+\ell)H$ for in $W_\epsilon^{j\ell}(t,\x,\K)$, and at $\xn= -x_{0,n} - (j+\ell)H$ for $W_\epsilon^{j\ell}(t,(\xperp,-x_n), (\Kperp,-k_n))$.
In the remainder of the paper, the terms $V_\epsilon^{j\ell}$ and $W_\epsilon^{j\ell}$ will be referred to as \emph{cross-Wigner transforms}. 

For further reference, the first sum in Eq.~\eqref{eq:wigner_total} will be denoted
\begin{multline}
 W_\epsilon^{1,\sharp}(t,\x,\K) = \sum_{j\in \mathbb{Z}} W_\epsilon\left(t,(\xperp,x_n-2jH), \K\right) \\+  W_\epsilon\left(t,(\xperp,-x_n+2jH),(\Kperp,-k_n) \right).
\end{multline}
As we will see below, in the high-frequency limit $\epsilon\to 0$, this sum will give rise to the energy density propagation over the domain, whose description will be given by a RTE similar to the one obtained in a full-space. The other two sums in Eq.~\eqref{eq:wigner_total} will be denoted
\begin{equation}
 W_\epsilon^{2,\sharp}(t,\x,\K) = \sum_{j,\ell \in \mathbb{Z} } V_\epsilon^{j\ell}(t,\x,\K) + V_\epsilon^{j\ell}\left(t,(\xperp,-x_n), (\Kperp,-k_n)\right),
\end{equation}
and
\begin{equation}
 W_\epsilon^{3,\sharp}(t,\x,\K) = \sum_{\substack{ j,l \in \mathbb{Z}   \\ \ell \ne j}} W_\epsilon^{j\ell}(t,\x,\K) + W_\epsilon^{j\ell}\left(t,(\xperp,-x_n), (\Kperp,-k_n)\right).
\end{equation}
The cross terms are responsible for interference effects between two waves arising from different initial conditions. More specifically, the second sum $W_\epsilon^{2,\sharp}$ will give rise to amplifications of the total energy in the vicinity (within one wavelength) of the boundaries, while the third sum $W_\epsilon^{3,\sharp}$ will lead to weak-localization effects in the vicinity of the planes $\left\{x_n = x_{0,n}\right\}$ and $\left\{x_n = H-x_{0,n}\right\}$.



\section{Asymptotic analysis of the Wigner transform}
\label{sec:asymptotic_analysis_with_an_initial_condition_far_from_the_border}
In this section, we examine the asymptotic limit of the three terms in Eq.~\eqref{eq:wigner_total}. We first show that the energy density propagation in the slab can be described through the limit of $W_\epsilon^{1,\sharp}$ (restricted to $\Omega$). Second, we exhibit interference effects at the boundaries and within the slab through an asymptotic analysis of the cross terms of Eq.~\eqref{eq:wigner_total}. 


\subsection{Radiative transfer equations for the self-Wigner transforms}
\label{sec:RTE_for_the_self_Wigner_transforms}

As $W_\epsilon$ corresponds to a problem posed on the entire space, the derivation of the radiative transfer model for $W_\epsilon$ is based on a classical multi-scale asymptotic analysis that can be found for instance in~\cite{Bal2005, butz2015}. The only difference in the present situation is the lack of stationarity of the random field modeling the mechanical properties. However, as in~\cite{Messaoudi2022}, where the mechanical properties were not periodic but even with respect to $\xn=0$, this lack of stationarity has no leading-order influence and the asymptotic analysis for the self-Wigner transforms $W_\epsilon$ remains similar to the classical case (for instance~\cite{Bal2005}). Thus, the limit self-Wigner transform can be written:
\begin{equation}
\label{eq:W_0_self}
W_0(t,\x,\K) := \lim_{\epsilon \to 0} W_\epsilon(t,\x,\K)
= a (t,\x,\K)\BB(\K) + a (t,\x,-\K)\BB^T(\K),
\end{equation}
where 
\begin{equation}\label{def:B}
\BB(\K):= \frac{1}{2}\begin{pmatrix}
1/\vert \K\vert^2 & i/(c_0\vert \K\vert) \\
- i/(c_0\vert \K\vert) & 1/c_0^2
\end{pmatrix},
\end{equation}
and $\BB^T$ stands for the transposition of $\BB$.
The decomposition~\eqref{eq:W_0_self} is an expansion over a basis of projectors associated to the eigenvectors of the dispersion matrix \[
\Pdiff_0(\I \K) = - \begin{pmatrix}
0 & c_0^2 \\
-\vert \K\vert^2 & 0
\end{pmatrix}
\]
of $\mathcal{A}_\epsilon$ in the hyperbolic system~\eqref{eq:equation_hyperbolique_ordre_1}. This decomposition is derived through the analysis presented in detail in Appendix \ref{sec:appendix_spectral_decomposition}.
The amplitude $a(t,\x,\K)$ represents the energy density carried by the wave field $q_\epsilon$ at time $t$, position $\x$ and direction $\K$ in the high-frequency limit $\epsilon\to 0$. This energy density evolves according to the following RTE 
\begin{multline}
\label{eq:equation_de_transport}
\partial_t  a (t,\x,\K)  + c_0 \widehat{\K} \cdot \nabla_\x a (t,\x,\K)   = -\Sigma(\K) a (t,\x,\K)  
\\+\int_{\mathbb{R}^3}\sigma(\K,\q) a (t,\x,\q) \delta(c_0(\vert \q \vert -\vert\K\vert ))d\q,
\end{multline}
equipped with the initial condition 
\begin{equation}
\label{eq:a0}
a(t=0,\x,\K) = \mathbb{A}(\K) \delta(\x-\x_{0}),
\end{equation}
where 
\begin{equation}
\label{def:A}
\mathbb{A}(\K) := \frac{1}{2(2\pi)^3}\Big\vert c_0^{-1} \widehat{B}(\K)  -i \vert \K \vert \widehat{A}(\K) \Big\vert^2.
\end{equation}  
In Eq.~\eqref{eq:equation_de_transport} 
\begin{equation}
\label{eq:scattering_coefficiant_1}
\Sigma(\K):=\frac{\pi c_0^2\vert \K \vert^2}{2(2\pi)^3}\int_{\mathbb{R}^3} \widehat{R}(\K-\q) \delta(c_0(\vert \q \vert -\vert\K\vert ))d\q, 
\end{equation}
and 
\begin{equation}
\sigma(\K,\q):=\frac{\pi c_0^2\vert \K \vert^2}{2(2\pi)^3}\widehat{R}(\K-\q).
\end{equation}
As a result, the overall contribution $W_\epsilon^{1,\sharp}$, in Eq.~\eqref{eq:wigner_total}, of the self-Wigner transform is given by 
\begin{align}\label{eq:limW1}
\lim_{\epsilon\to 0} W_\epsilon^{1,\sharp} (t,\x,\K) & =\sum_{j\in \mathbb{Z}} \lim_{\epsilon\to 0} W_\epsilon\left(t,(\xperp, x_n -2jH  ),\K \right) \\
& \hspace{2cm}+ \lim_{\epsilon\to 0} W_\epsilon \left(t,(\xperp,-x_n +2jH),(\Kperp,-k_n) \right)\nonumber\\
& =\sum_{j\in \mathbb{Z}} W_0\left(t,(\xperp, x_n -2jH  ),\K \right) \\
& \hspace{2cm}+  W_0 \left(t,(\xperp,-x_n +2jH),(\Kperp,-k_n) \right).
\end{align}

Although explicit solutions to Eq.~\eqref{eq:equation_de_transport} can be exhibited for very particular problems (see \cite{Paasschens1997} for instance), solutions in the general case must be obtained numerically, for instance using finite difference or finite elements, Lattice Boltzmann Methods~\cite{Mink2020,Wang2021,Ruyssen2023}, neural networks~\cite{Mishra2021}, or, more adapted to the dimensionality of the RTE equations, Monte-Carlo methods~\cite{Lapeyre2003} (see also the github repository \texttt{https://github.com/cottereau/RadiativeTransferMonteCarlo}, that implements the method associated to the configurations described in~\cite{Messaoudi2022} and in the present paper).

Following the same lines as in \cite{Messaoudi2022} regarding the half-plane, one can show that, without any particular precaution, the cross-Wigner transforms $V_\epsilon^{j\ell}$ and $W_\epsilon^{j\ell}$ 
go to $0$ as $\epsilon\to 0$. In fact, both of the limit of $V_\epsilon^{j\ell}$ and $W_\epsilon^{j\ell}$ can be described as verifying standard (linear) RTEs with vanishing initial conditions. As a result, taking the high-frequency limit $\epsilon\to 0$ in Eq.~\eqref{eq:wigner_total}, the total energy density is only given by
\begin{multline}
W^{tot} (t,\x,\K) := \lim_{\epsilon\to 0} W_\epsilon^\sharp (t,\x,\K) =\lim_{\epsilon\to 0} W_\epsilon^{1,\sharp} (t,\x,\K) \\
 =\sum_{j\in \mathbb{Z}} W_0\left(t,(\xperp, x_n -2jH  ),\K \right) +  W_0 \left(t,(\xperp,-x_n +2jH),(\Kperp,-k_n) \right).
\end{multline}
From Eq.~\eqref{eq:W_0_self}, the full energy density $W^{tot}$ can be expressed as 
\begin{align}
\label{eq:Wtot}
W^{tot} (t,\x,\K) &= \sum_{j\in \mathbb{Z}} \left( a\left( t,(\xperp,x_n-2jH),\K \right) \right. \nonumber\\
& \hspace{2cm} \left.+ a\left( t,(\xperp,-x_n+2jH),(\Kperp,-k_n) \right) \right)\BB(\K) \nonumber\\
&+ \sum_{j\in \mathbb{Z}} \left( a\left(t,(\xperp,x_n-2jH),-\K \right) \right.\\
& \hspace{2cm} \left.+ a\left( t,(\xperp,-x_n+2jH),(-\Kperp,k_n) \right) \right)\BB^T(\K). \nonumber
\end{align}
Using arguments similar to those of Sec.~\ref{sec:the_method_of_images}, this superposition corresponds to imposing reflecting boundary conditions~\eqref{eq:reflecting_boundary_conditions} on the field $W^{tot}(t,(\xperp, x_n))$ at the boundary $\{x_n=0\}$:
\begin{equation}\label{eq:reflecting_boundary_conditions0}
W^{tot} (t,(\xperp,x_n=0),\K) = W^{tot} (t,(\xperp,x_n=0),(\Kperp,-k_n)).
\end{equation}
This boundary condition corresponds to the one obtained in \cite{Jin2006b,Gerard1993, Ryzhik1997, Burq1997}, and a similar condition holds at the interface $\{x_n=H\}$
\begin{equation}\label{eq:reflecting_boundary_conditionsH}
W^{tot} (t,(\xperp,x_n=H),\K) = W^{tot} (t,(\xperp,x_n=H),(\Kperp,-k_n)).
\end{equation}

Finally, introducing the total energy contributions from all directions
\begin{equation}
\label{eq:energy1}
E(t,\x) := \int_{\mathbb{R}^3} W^{tot} (t,\x,\K) d\K, 
\end{equation}
we have
\begin{align*}
E(t,\x) & = \sum_{j\in \mathbb{Z}} \int_{\mathbb{R}^3} W_0\left(t,(\xperp, x_n -2jH  ),\K \right) \\
& \hspace{2cm}+ W_0\left(t,(\xperp,-x_n +2jH),(\Kperp,-k_n) \right) d\K \\
& =\sum_{j\in \mathbb{Z}} \int_{\mathbb{R}^3} \left( a\left( t,(\xperp,x_n-2jH),\K \right) \right.\\
& \hspace{2cm} \left.+ a\left( t,(\xperp,-x_n+2jH),\K \right) \right) \DD(\K) d\K,
\end{align*}
where
\begin{equation}
\label{def:D}
\DD(\K):= \begin{pmatrix}
1/\vert \K\vert^2 & 0 \\
0 & 1/c_0^2
\end{pmatrix},
\end{equation}
knowing that $\BB(\K) = b(\vert\K\vert)$.

To summarize this section, it has been shown that, following the standard asymptotic analysis, only the self-Wigner transforms contribute to the total energy over the whole slab $\Omega$. At leading order in $\epsilon$, no nontrivial contributions of the cross-terms have been observed. 

\subsection{Energy contribution of the cross-terms}
\label{sec:energy_contribution_of_the_cross_terms}

To complete the analysis of the full energy density $W_\epsilon^\sharp$, it remains to describe the asymptotic behavior of the other terms in Eq.~\eqref{eq:wigner_total}: the cross-Wigner transforms in $W_\epsilon^{2,\sharp}$ and $W_\epsilon^{3,\sharp}$. To highlight the contributions provided by these terms and exhibit the intensity enhancement and the weak-localization phenomena, we need a more careful approach compared to the standard asymptotic analysis.


\subsubsection{Asymptotic behavior of the cross-Wigner transforms $V_\epsilon^{j \ell}$}
\label{sec:contribution_v_eps}
In the following, for notational simplicity, the slow component in $c_\sharp$ will be implicit, so that by $c_\sharp(\x/\epsilon)$ we mean $c_\sharp(\x,\x/\epsilon)$. To exhibit the nontrivial contributions carried by the cross-Wigner transforms $V_\epsilon^{j \ell}$, we need to take a close look at these quantities. It turns out that $V_\epsilon^{j \ell}$ can be rewritten as 
\begin{multline}
V_\epsilon^{j \ell}(t,\x,\K)  = \int_{\mathbb{R}^3}\e^{\I\K\cdot \y}  \U_\epsilon\big(t,\x_\perp-\epsilon\y_\perp/2,  -2jH   - \epsilon y_n/2  + (j+\ell)H   + ( x_n - (j+\ell)H)   \big) \\
 \otimes \U_\epsilon \big(t,\x_\perp + \epsilon\y_\perp/2,   -2jH   - \epsilon y_n/2  + (j+\ell)H  - ( x_n - (j+\ell)H) \big) \frac{d\y}{(2\pi) ^3}.
\end{multline}
In this way, the above cross-Wigner transform can be recast as
\begin{align}\label{eq:cross_Wigner_1}
V_\epsilon^{j \ell}(t,\x,\K) & = \frac{1}{2\pi} \int_{\mathbb{R}} \e^{\I k_ny_n} \Big( \int_{\mathbb{R}} \e^{2\I p_n ( x_n - (j+\ell)H )/\epsilon} \nonumber \\
&\hspace{2cm} \times W_\epsilon \left(t, ( \x_\perp,  -2jH  - \epsilon y_n/2 + (j+\ell)H   )   , (\Kperp,p_n) \right) dp_n \Big)dy_n\nonumber\\
& = \frac{2}{\epsilon(2\pi)} \int_{\mathbb{R}} \int_{\mathbb{R}} \e^{2 \I k_n(y_n- 2jH)/\epsilon}  \e^{2\I p_n  ( x_n - (j+\ell)H)/\epsilon} \nonumber\\
&\hspace{2cm} \times W_\epsilon\left(t, ( \xperp,  (j+\ell)H - y_n )   , (\Kperp,p_n) \right) dp_ndy_n,
\end{align}
\noindent where $W_\epsilon$ is defined by Eq.~\eqref{eq:self_wigner}.
From this latter relation let us make some remarks. First, although the asymptotic ($\epsilon\to 0$) initial condition of $V_\epsilon^{j\ell}$ is $0$, this is no more the case for $W_\epsilon$. 
As a result, the standard asymptotic analysis for this term will provide a nontrivial contribution. Second, the presence of highly oscillatory terms $\e^{2\I p_n  ( x_n - (j+\ell)H)/\epsilon}$ suggests to place ourselves in the vicinity of the planes $\left\{ \xn = (j+\ell)H)  \right\}$ (within one wavelength), leading to the changes of variables 
\begin{equation}\label{def:chg_var1}
\x \to \x_\epsilon^{j\ell}:= \big(\xperp,(j+\ell)H + \epsilon \tilde{x}_n/2 \big).
\end{equation}
This change of variable allows us to focus on interference phenomena near the planes $\{ \xn = (j+\ell)H \}$. Although the propagation medium considered throughout our analysis is a 3D full-space, the objective is to describe the intensity enhancement phenomena in the slab $\mathbb{R}^2\times[0,H]$. As a result, we necessarily have 
\[
\ell=-j \qquad \text{and} \qquad \ell=1-j,
\]
corresponding to the interface $\{x_n=0\}$ and $\{x_n=H\}$ respectively. All the other choices for $\ell$ would lead to negligible contribution due to a remaining integrated fast phase that average out in the limit $\epsilon\to 0$.  
Third, the oscillatory term $\e^{2\I k_n(y_n-(x_{0,n} + 2jH))/\epsilon}$ leads to the change of variable 
\begin{equation}
\label{eq:change_of_var}
\K \to  \K_\epsilon :=\big(\Kperp,\epsilon \, \tilde{k}_n/2\big),
\end{equation}
meaning that we focus on energy propagation along the planes $\{ \xn = (j+\ell)H \}$. The above changes of variables yields an asymptotic nontrivial contribution
\begin{multline}
\label{eq:asymptotic_limit_cross_wigner_1}
 \lim \limits_{\epsilon\to 0} \epsilon V_\epsilon^{j \ell} (t, \x_\epsilon^{j\ell},\K_\epsilon) 
  =  \frac{2}{(2\pi)} \int_{\mathbb{R}} \int_{\mathbb{R}} \e^{\I \tilde{k}_n(y_n-2jH )}  \e^{ \I p_n \tilde{x}_n}\\ \times W_0 \left(t,  \xperp, (j+\ell)H -y_n  , (\Kperp,p_n )\right)  dp_ndy_n.    
\end{multline}
where $W_0$ is given by Eq.~\eqref{eq:W_0_self}.

The same analysis as above leads to the following expression for $V_\epsilon^{j \ell}(t,(\xperp,-x_n),(\Kperp,-k_n))$ 
\begin{multline*}
V_\epsilon^{j \ell}(t,(\xperp,-x_n),(\Kperp,-k_n)) 
=  \frac{2}{\epsilon(2\pi)} \int_{\mathbb{R}} \int_{\mathbb{R}} \e^{-2\I k_n(y_n-  2jH)/\epsilon}  \e^{2\I p_n  ( -x_n - (j+\ell)H)/\epsilon} \\
\times W_\epsilon \left(t, ( \xperp, (j+\ell)H - y_n )   , (\Kperp,p_n) \right) dp_ndy_n,       
\end{multline*} 
suggesting the changes of variables
\begin{equation*}
\x \to \tilde{\x}_\epsilon^{j\ell}:=\left(\xperp,-(j+\ell)H + \epsilon \, \tilde{x}_n/2    \right) \quad \mathrm{and} \quad \K \to  \K_\epsilon.
\end{equation*} 
This time, the change of variable $\x \to \tilde{\x}_\epsilon^{j\ell}$ allows us to focus on interference phenomena near the planes $\left\{ \xn = -(j+\ell)H  \right\}$, but to restrict our attention to the slab $\mathbb{R}^2\times [0,H]$, we necessarily have 
\[
\ell=-j \qquad \text{and} \qquad \ell=-1-j,
\]
corresponding to the interface $\{x_n=0\}$ and $\{x_n=H\}$ respectively. In the same way as Eq.~\eqref{eq:asymptotic_limit_cross_wigner_1}, we deduce that
\begin{multline}
\label{eq:asymptotic_limit_cross_wigner_1_bis}
\lim \limits_{\epsilon\to 0} \epsilon V_\epsilon^{j \ell} (t, (\xperp,(j+\ell)H - \epsilon \tilde{x}_n/2),(\Kperp,-\epsilon \tilde k_n/2)) \\
=  \frac{2}{(2\pi)} \int_{\mathbb{R}} \int_{\mathbb{R}} \e^{-\I \tilde{k}_n(y_n-2jH )}  \e^{- \I p_n \tilde{x}_n} 
 W_0 \left(t,  \xperp, (j+\ell)H -y_n , (\Kperp,p_n )\right) dp_ndy_n.    
\end{multline}
As described in the next section, this analysis is the key tool to highlight the intensity enhancement phenomena occurring within one wavelength of the slab boundaries.

\subsubsection{Intensity enhancement at the boundaries}
\label{sec:amplification_at_the_booundaries_for_the_cross_Wigner_transforms}

Focusing firstly on the boundary $\{x_n = 0\}$, which necessarily yields the relation $\ell = -j$ as pointed out in the previous section, $\x_\epsilon^{j\ell}$ and $\tilde{\x}_\epsilon^{j\ell}$ just become $\x_\epsilon = (\x_{\perp},\epsilon \tilde x_n /2)$. The energy contributions within one wavelength of the boundary then read
\begin{equation*}
E_{\{x_n = 0\}}(t,\xperp,\tilde {x}_n) := \int_{\mathbb{R}^3} \lim_{\epsilon\to 0} W_\epsilon^\sharp(t,\x_\epsilon,\K) d\K
=\int_{\mathbb{R}^3} \lim_{\epsilon\to 0} W_\epsilon^\sharp(t,(\xperp,\epsilon \tilde x_n/2),\K) d\K.
\end{equation*}
From Eq.~\eqref{eq:wigner_total}, together with the asymptotic analysis for $V_\epsilon^{j \ell}$, we obtain, after the change of variable $\K\to \K_\epsilon$ for the integral of the cross-terms, 
\begin{align*}
E_{\{x_n = 0\}}(t,\xperp,\tilde {x}_n) & = \sum_{j\in\mathbb{Z}} \int_{\mathbb{R}^3} \lim_{\epsilon\to 0} W_\epsilon(t,(\xperp, -2jH+\epsilon \tilde x_n/2 ),\K) \\
& \hspace{1.3cm}+ \lim_{\epsilon\to 0}W_\epsilon(t,(\xperp,2jH-\epsilon \tilde x_n/2 ),(\Kperp,-k_n)) d\K \\
& +  \frac{1}{2}  \sum_{ j\in \mathbb{Z} } \int_{\mathbb{R}^3} \lim_{\epsilon\to 0} \epsilon  V_\epsilon^{j, \ell = -j} (t, (\xperp,\epsilon \tilde x_n/2),(\Kperp, \epsilon \tilde k_n /2)) \\
& \hspace{1.3cm} + \lim_{\epsilon\to 0} \epsilon V_\epsilon^{j, \ell=-j} (t, (\xperp,-\epsilon \tilde x_n/2),(\Kperp,-\epsilon \tilde k_n /2)) d\K_\perp d\tilde k_n.
\end{align*}
In this latter relation, one can see that the self-Wigner transforms are not the only terms contributing to the limit. Note that, as we will see in the next section, the contribution of $W_\epsilon^{3,\sharp}$ (in Eq.~\eqref{eq:wigner_total}) goes to $0$ in the limit $\epsilon\to 0$ when we look at the energy close to the slab boundaries. This term produces a nontrivial contribution only within the slab. Now, remembering Eq.~\eqref{eq:W_0_self}, Eq.~\eqref{eq:asymptotic_limit_cross_wigner_1} and Eq.~\eqref{eq:asymptotic_limit_cross_wigner_1_bis}, the total energy becomes
\begin{align*}
E_{\{x_n = 0\}}(t,\xperp,\tilde {x}_n) &=\sum_{j\in\mathbb{Z}} \int_{\mathbb{R}^3} W_0(t,(\xperp, -2jH),\K) \\
&\hspace{2cm} + W_0(t,(\xperp,2jH),(\Kperp,-k_n)) d\K \\
&\hspace{0.2cm}+ \sum_{j\in \mathbb{Z}}  \int_{\mathbb{R}^3} \e^{ \I k_n \tilde{x}_n} W_0 \left(t,  (\xperp, -2jH )  , \K \right)  d\K \\
&\hspace{2cm} +\sum_{j\in \mathbb{Z}}  \int_{\mathbb{R}^3} \e^{- \I k_n \tilde{x}_n} W_0 \left(t,  (\xperp, -2jH)   , \K\right)  d\K,
\end{align*}
which can be recast as
\begin{equation}
\label{eq:energy_contribution_0}
E_{\{x_n = 0\}}(t,\xperp,\tilde {x}_n) = 2 \sum_{j\in \mathbb{Z}}  \int_{\mathbb{R}^3}a\left( t,(\xperp,2jH),\K \right)  (1+\cos(k_n \tilde{x}_n))  \DD(\K) d\K, 
\end{equation}
according to Eq.~\eqref{eq:W_0_self}, where $\DD$ is given by Eq.~\eqref{def:D}. At exactly the boundary, that is for $\tilde {x}_n = 0$, we therefore observe a doubling of the total energy w.r.t. Eq.~\eqref{eq:energy1}:
\[E_{\{x_n = 0\}}(t,\xperp,\tilde {x}_n = 0) = 2 \, E(t,(\xperp,0)).\]

Regarding the boundary $\{x_n = H\}$, we proceed in a very similar way, but this time with
\begin{align*}
E_{\{x_n = H\}}(t,\xperp,\tilde {x}_n) & = \sum_{j\in\mathbb{Z}} \int_{\mathbb{R}^3} \lim_{\epsilon\to 0} W_\epsilon(t,(\xperp, H-2jH+\epsilon \tilde x_n/2 ),\K) \\
& \hspace{1.5cm}+ \lim_{\epsilon\to 0}W_\epsilon(t,(\xperp,-H + 2jH-\epsilon \tilde x_n/2 ),(\Kperp,-k_n)) d\K \\
& \hspace{0.1cm}  +  \frac{1}{2}  \sum_{ j\in \mathbb{Z} } \int_{\mathbb{R}^3} \lim_{\epsilon\to 0} \epsilon  V_\epsilon^{j, \ell = 1-j} (t, (\xperp,H+\epsilon \tilde x_n/2),(\Kperp,\epsilon \tilde k_n /2)) \\
& \hspace{0.2cm} + \lim_{\epsilon\to 0} \epsilon V_\epsilon^{j, \ell=-1-j} (t, (\xperp,H-\epsilon \tilde x_n/2),(\Kperp,-\epsilon \tilde k_n /2)) d\K_\perp d\tilde k_n,
\end{align*}
yielding 
\begin{equation}\label{eq:energy_contribution_H}
E_{\{x_n = H\}}(t,\xperp,\tilde {x}_n) = 2 \sum_{j\in \mathbb{Z}}  \int_{\mathbb{R}^3}  a\left( t,(\xperp,(2j+1)H),\K \right) (1+\cos(k_n \tilde{x}_n))  \DD(\K) d\K.
\end{equation}
Again, at exactly the boundary (that is for $\tilde {x}_n = 0$), we observe a doubling of the total energy w.r.t. Eq.~\eqref{eq:energy1},
\[E_{\{x_n = H\}}(t,\xperp,\tilde {x}_n = 0) = 2 \, E(t,(\xperp,H)).\]



\subsubsection{Asymptotic behavior of the cross-Wigner transforms $W_\epsilon^{j \ell}$}
\label{sec:contribution_W_eps}
 
To exhibit the nontrivial contributions provided by the cross-Wigner transforms $W_\epsilon^{j \ell}$ we follow the same lines as in Sec.~\ref{sec:contribution_v_eps}. The first step is to rewrite $W_\epsilon^{j \ell}$ in a more convenient form using the \emph{shifted} wave field 
\begin{equation}
\label{eq:translation_scalar_wave_fields}
p^{j \ell}_\epsilon(t,\x) = q_\epsilon(t,\xperp, \xn - (j+\ell)H ).
\end{equation}
In other words, $p^{j \ell}_\epsilon$ represents the wave propagating in the periodically-structured full-space, in the weak-scattering regime whose initial condition is centered around $(\x_{0,\perp},x_{0,n}+(j+\ell)H)$. We also consider the \emph{shifted} wave field $g^{j \ell}_\epsilon$ whose initial condition is the same as that of $p^{j \ell}_\epsilon$. However, the propagation media of these two wave fields are shifted by $2x_{0,n}$ w.r.t the $x_n$-variable. More precisely, $g^{j \ell}_\epsilon$ satisfies the following \emph{shifted} wave equation
\begin{equation}
\label{eq:equation_onde_translate_1}
\partial^2_{tt} g^{j\ell}_\epsilon(t,\x)-\ctilde^2\left(\frac{\x_\perp}{\epsilon},\frac{x_n -2x_{0,n} - (j+\ell)H}{\epsilon}\right) \Delta g^{j\ell}_\epsilon(t,\x)  = 0 \qquad (t,\x)\in\mathbb{R}_+^\ast \times \mathbb{R}^3,
\end{equation} 
equipped with the initial conditions
\[
g_\epsilon^{j\ell}(t=0,\x)= p_\epsilon^{j\ell}(t=0,\x)=\frac{1}{\epsilon^{3/2}}A\left(  \frac{\x_\perp- \x_{0,\perp}}{\epsilon}, \frac{x_n-x_{0,n}-(j+\ell)H}{\epsilon}   \right)\]
and
\[
\epsilon \partial_t g_\epsilon^{j\ell}(t=0,\x)=\epsilon \partial_t p_\epsilon^{j\ell}(t=0,\x)=\frac{1}{\epsilon^{3/2}}B\left(  \frac{\x_\perp- \x_{0,\perp} }{\epsilon},  \frac{x_n-x_{0,n}-(j+\ell)H}{\epsilon}   \right).
\]

Considering now the two vectors
\begin{align}
\label{eq:cross_waves_vectors_axis}
\mathbf{g}_\epsilon^{j \ell}(t,\x)&:=\begin{pmatrix} g_\epsilon^{j\ell}(t,\x) \\ \epsilon\ctilde^{-2}\left(\frac{\x_\perp}{\epsilon},\frac{x_n - 2x_{0,n} - (j+\ell)H}{\epsilon}\right) \partial_t g_\epsilon^{j\ell}(t,\x)  \end{pmatrix}, \nonumber \\
&\\
 \p_\epsilon^{j\ell}(t,\x)& :=  \U_\epsilon(t,\xperp,x_n-(j+\ell) H),\nonumber  
\end{align} 
and remembering that $A$ and $B$ are assumed to be even w.r.t the $\xn$-variable for simplicity, as well as $\ctilde$ is $2H$-periodic and even w.r.t the $\xn$-variable, we deduce that
\[
\U_\epsilon(t,\x_\perp,x_n-2jH)=\mathbf{g}_\epsilon^{j\ell}\big(t,\x_\perp,-x_n + (x_{0,n} + 2jH) + ( x_{0,n} + (j + \ell)H) \big), 
\]
and
\[
\U_\epsilon(t,\x_\perp,x_n-2\ell H)=\p_\epsilon^{j\ell}\big(t,\x_\perp,x_n  + (x_{0,n} + 2 j H) - ( x_{0,n} + (j + \ell)H) \big).
\]
The cross-Wigner transform $W_\epsilon^{j \ell}$ can then be rewritten as
\begin{align*}
W_\epsilon^{j \ell}(t,\x,\K) & = \int_{\mathbb{R}^3}\e^{\I\K\cdot \y}  \mathbf{g}_\epsilon^{j\ell}\big(t,\x_\perp-\epsilon\y_\perp/2,-x_n+\epsilon y_n/2 + (x_{0,n} + 2jH) + ( x_{0,n} + (j + \ell)H)  \big)  \\
& \hspace{1.5cm}\otimes \p_\epsilon^{j\ell} \big(t,\x_\perp + \epsilon\y_\perp/2,  x_n+\epsilon y_n/2 + (x_{0,n} + 2 j H) - ( x_{0,n} + (j + \ell)H)    \big) \frac{d\y}{(2\pi) ^3} \\
& = \int_{\mathbb{R}^3}\e^{\I\K\cdot \y}   \mathbf{g}_\epsilon^{j\ell}\big(t,\x_\perp-\epsilon\y_\perp/2, x_{0,n} + 2jH   + \epsilon y_n/2  -( x_n - (x_{0,n}+ (j+\ell)H))  \big) \\ 
& \hspace{1.5cm}\otimes \p_\epsilon^{j \ell} \big(t,\x_\perp + \epsilon\y_\perp/2,  x_{0,n} + 2jH   + \epsilon y_n/2  + ( x_n - \big(x_{0,n}+ (j+\ell)H )) \big) \frac{d\y}{(2\pi) ^3}.
\end{align*}
In the same way as Eq.~\eqref{eq:cross_Wigner_1}, $W_\epsilon^{j \ell}$ can be expressed as follows
\begin{multline*}
W_\epsilon^{j \ell}(t,\x,\K) =  \frac{2}{\epsilon(2\pi)} \int_{\mathbb{R}} \int_{\mathbb{R}} \e^{2\I k_n(y_n-(x_{0,n} + 2jH))/\epsilon}  \e^{-2\I p_n  ( x_n - (x_{0,n}+ (j+\ell)H))/\epsilon} \\
\times \W_\epsilon^{j\ell}\left(t, ( \xperp, y_n )   , (\Kperp,p_n) \right) dp_ndy_n,      
\end{multline*} 
in terms of the cross-Wigner transform
\begin{equation}
\label{def:calW}
\W_\epsilon^{j\ell}(t,\x,\K) = W[\mathbf{g}_\epsilon^{j\ell},\p_\epsilon^{j\ell}](t,\x,\K). 
\end{equation}
As for $V_\epsilon^{j\ell}$, despite $W_\epsilon^{j\ell}$ has trivial asymptotic initial conditions, this is no more the case for $\mathcal{W}_\epsilon^{j\ell}$.  The standard asymptotic analysis therefore provides a nontrivial contribution for this term. Moreover, the presence of the highly oscillatory terms $\e^{-2\I p_n  ( x_n - (x_{0,n}+ (j+\ell)H))/\epsilon}$ in the latter representation of $W_\epsilon^{j \ell}$ leads to the change of variables  
\begin{equation*}
\x \to \z_\epsilon^{j\ell}:=\big(\xperp, x_{0,n} + (j+\ell)H + \epsilon \, \tilde{x}_n/2    \big),
\end{equation*} 
focusing our attention around the planes passing through $x_n = x_{0,n} + (j+\ell)H$. To highlight the contribution of $W_\epsilon^{j, \ell}$ within the slab $\mathbb{R}^2\times [0,H]$, we necessarily have
\[
\ell = -j.
\]
All other choices of $\ell$ lead to a negligible contribution. Note that this also excludes the case $j=0$, since we would have $\ell=0$, which is in contradiction with $\ell\neq j$ in Eq.~\eqref{eq:wigner_total}. The other rapid phase $\e^{2\I k_n(y_n-(x_{0,n} + 2jH))/\epsilon}$ leads again to the change of variable $\K \to \K_\epsilon$, where $\K_\epsilon$ is given by Eq.~\eqref{eq:change_of_var}, meaning that we look at the energy propagating along the plane $\{x_n=x_{0,n}\}$ (within one wavelength). As a result, we obtain
\begin{multline*}
 W_\epsilon^{j, \ell=-j} (t, \z_\epsilon^{j,\ell=-j},\K_\epsilon) 
=  
\frac{2}{\epsilon(2\pi)} \int_{\mathbb{R}} \int_{\mathbb{R}} \e^{\I \tilde{k}_n(y_n-(x_{0,n}+2jH) )}  \e^{- \I p_n \tilde{x}_n} \\\times \W_\epsilon^{j,\ell=-j} \left(t, ( \xperp, y_n )   , (\Kperp,p_n )\right)  dp_ndy_n,   
\end{multline*}
and the asymptotic (nontrivial) contribution of $W_\epsilon^{j \ell} $ is given by
\begin{multline}
\label{eq:asymptotic_limit_cross_wigner_2}
 \lim \limits_{\epsilon\to 0} \epsilon W_\epsilon^{j ,\ell=-j} (t, \z_\epsilon^{j,\ell=-j},\K_\epsilon) 
=  
\frac{2}{(2\pi)} \int_{\mathbb{R}} \int_{\mathbb{R}} \e^{\I \tilde{k}_n(y_n-(x_{0,n}+2jH) )}  \e^{- \I p_n \tilde{x}_n}\\
\times \W_0^{j,\ell=-j} \left(t, ( \xperp, y_n )   , (\Kperp,p_n )\right)  dp_ndy_n.    
\end{multline}
Here, we have
\begin{equation}
\label{def:W_0_jl}
\W_0^{j\ell}(t,\x,\K) :=  \lim \limits_{\epsilon\to 0} \W_\epsilon^{j\ell}(t,\x,\K) = a^{j\ell}(t,\x, \K) \BB(\K)  + a^{j\ell}(t,\x,-\K)\BB^T(\K),
\end{equation}
where $\BB$ is given by Eq.~\eqref{def:B}, and $a^{j\ell}$ satisfies the following transport equation 
\begin{equation}
\label{eq:RTE_coherent}
\partial_t  a ^{j\ell} (t,\x,\K)  + c_0 \widehat{\K} \cdot \nabla_\x a ^{j\ell} (t,\x,\K)   = -\Sigma(\K) a  ^{j\ell} (t,\x,\K),
\end{equation}
equipped with the initial condition 
\begin{equation}
\label{eq:a0_jl}
a^{j\ell}(t=0,\x,\K) = a^{j\ell}_0(\x,\K) := \mathbb{A}(\K) \delta(\xperp-\x_{0,\perp})\delta(x_n-x_{0,n}-(j+\ell)H).
\end{equation}
$\mathbb{A}$ is given by Eq.~\eqref{def:A} and $\Sigma$ is given by Eq.~\eqref{eq:scattering_coefficiant_1}.
Such an equation admits an explicit solution given by
\begin{equation}\label{eq:explicite_form_ajl}
a^{j\ell}(t,\x,\K) = \e^{-\Sigma(\K)t}a^{j\ell}_0(\x-c_0t\widehat{\K},\K).  
\end{equation}
The derivation of the transport equation~\eqref{eq:RTE_coherent} from $\W_\epsilon^{j\ell}$ is detailed in Appendix \ref{sec:appendix}. This derivation holds actually for the average quantities $\langle \W_\epsilon^{j\ell}\rangle$, but a deeper analysis such as in \cite{butz2015} shows that the (deterministic) transport limit holds also for $\W_\epsilon^{j\ell}$ itself. This \emph{self-averaging} property, appearing also for Eq.~\eqref{eq:equation_de_transport}, is due to the weak-scattering regime we consider in this paper. Note that compared to Eq.~\eqref{eq:equation_de_transport}, we do not have a transfer term anymore (the integral term), only a dissipation term. The reason is that, through Eq.~\eqref{def:calW}, we correlate two wave functions that have not seen the same propagation medium. The two propagation media of $g_\epsilon^{j\ell}$ and $p_\epsilon^{j\ell}$
are shifted by $2x_{0,n}$ w.r.t the $x_n$-variable. This offset being too large, the coda waves are too different to provide effective interferences through correlations. The resulting contribution Eq.~\eqref{eq:explicite_form_ajl} corresponds only to the coherent part of the radiative transfer model~\eqref{eq:equation_de_transport}.  
These observations are consistent with the results of \cite{Garnier2013, Bal2005} in which cross-correlations of wave fields with shifted propagation medium are analyzed. 

Regarding the term $W_\epsilon^{j \ell}(t,(\xperp,-x_n),(\Kperp,-k_n))$, the same analysis as above leads to 
\begin{multline*}
W_\epsilon^{j \ell}(t,(\xperp,-x_n),(\Kperp,-k_n)) =   \frac{2}{\epsilon(2\pi)} \int_{\mathbb{R}} \int_{\mathbb{R}} \e^{-2\I k_n(y_n-(x_{0,n} + 2jH))/\epsilon}  \\
\times \e^{2\I p_n  (x_n +x_{0,n} + (j+\ell)H)/\epsilon} \W_\epsilon^{j\ell}\left(t, ( \xperp, y_n )   , (\Kperp,p_n) \right) dp_ndy_n,       
\end{multline*} 
suggesting the changes of variables
\begin{equation*}
\x \to \tilde{\z}_\epsilon^{j\ell}:= \big(\xperp,-x_{0,n}-(j+\ell)H + \epsilon \, \tilde{x}_n/2 \big) \qquad \mathrm{and} \qquad \K \to  \K_\epsilon,
\end{equation*} 
to focus on energy propagation along the planes $\{x_n = -x_{0,n} -(j+l)H\}$. Here, to describe the interference phenomena inside the slab, we necessarily have
\[
\ell=-1-j,
\] 
so that the corresponding contribution takes place (within one wavelength) along the plane $\{x_n = H-x_{0,n}\}$. As for Eq.~\eqref{eq:asymptotic_limit_cross_wigner_2}, we deduce that
\begin{multline}
\label{eq:asymptotic_limit_cross_wigner_2_bis}
 \lim \limits_{\epsilon\to 0} \epsilon W_\epsilon^{j,\ell=-1-j} (t, (\x_\perp, x_{0,n}-H - \epsilon \tilde{x}_n/2  ),(\K_\perp,-\epsilon \tilde k_n /2)) 
\\ =  
\frac{2}{(2\pi)} \int_{\mathbb{R}} \int_{\mathbb{R}} \e^{-\I \tilde{k}_n(y_n-(x_{0,n}+2jH) )}  \e^{ \I p_n \tilde{x}_n} \W_0^{j\ell} \left(t, ( \xperp, y_n )   , (\Kperp,p_n )\right)  dp_ndy_n.    
\end{multline}

As discussed in the next section, the limits~\eqref{eq:asymptotic_limit_cross_wigner_2} and~\eqref{eq:asymptotic_limit_cross_wigner_2_bis} are the key quantities exhibiting interference effects within the slab (along the planes $\{x_n = x_{0,n}\}$ and $\{x_n = H-x_{0,n}\}$). These interferences are interpreted as weak localization effects as they result from wave fields which have not seen the propagation medium in the same way.



\subsubsection{Weak localization phenomenon inside the slab} 
\label{sec:weak_localization_phenomena_into_the_slab}

According to the above analysis, interference phenomena occurring along the planes $\{x_n = x_{0,n}\}$ and $\{x_n = H - x_{0,n}\}$ (the dashed lines in Figure~\ref{fig:superposition_ondes}) can be described. First, along the plane $\{x_n = x_{0,n}\}$, we necessarily have the relation $\ell = -j$, so that $\z_\epsilon^{j\ell}$ becomes $\z_\epsilon = (\x_{\perp},x_{0,n} + \epsilon x_n /2)$. The resulting energy contributions from all directions read
\begin{multline*}
E_{\{x_n = x_{0,n}\}}(t,\xperp,\tilde {x}_n) := \int_{\mathbb{R}^3} \lim_{\epsilon\to 0} W_\epsilon^\sharp(t,\z_\epsilon,\K) d\K\\
=\int_{\mathbb{R}^3} \lim_{\epsilon\to 0} W_\epsilon^\sharp(t,(\xperp,x_{0,n} + \epsilon \tilde x_n/2),\K) d\K.
\end{multline*}
From Eq.~\eqref{eq:wigner_total}, together with the asymptotic analysis of $W_\epsilon^{j \ell}$, we obtain, after the change of variable $\K\to \K_\epsilon$ for the integral of the cross-terms,
\begin{multline}
E_{\{x_n = x_{0,n}\}}(t,\xperp,\tilde {x}_n)  = \sum_{j\in\mathbb{Z}} \int_{\mathbb{R}^3} \lim_{\epsilon\to 0} W_\epsilon(t,(\xperp,  x_{0,n}-2jH+ \epsilon \tilde x_n/2),\K)  \\
\hspace{1.5cm}+ \lim_{\epsilon\to 0}W_\epsilon\left(t,\left(\xperp,- x_{0,n}+2jH-\epsilon \tilde x_n/2 \right),(\Kperp,-k_n)\right) d\K \\
 \hspace{0.5cm} +  \frac{1}{2}  \sum_{ j\in \mathbb{Z}^\ast } \int_{\mathbb{R}^3} \lim_{\epsilon\to 0} \epsilon  W_\epsilon^{j, \ell = -j} (t, (\xperp,x_{0,n}+\epsilon \tilde x_n/2),\K_\epsilon) d\K_\perp d\tilde k_n.
\end{multline}
Within one wavelength of the plane $\{x_n = x_{0,n}\}$, the energy is not provided by the self-Wigner transform only, the term $W_\epsilon^{3,\sharp}$ (in Eq.~\eqref{eq:wigner_total}) provides here additional contributions. Let us remind that the terms $W_\epsilon^{2,\sharp}$ goes to $0$ as $\epsilon\to 0$, since it provides non-negligible contributions only at the boundaries of the slab. Remembering Eq.~\eqref{eq:W_0_self}, Eq.~\eqref{eq:asymptotic_limit_cross_wigner_2}, and Eq.~\eqref{def:W_0_jl}, we deduce that
\begin{multline}
E_{\{x_n = x_{0,n}\}}(t,\xperp,\tilde {x}_n)  = \sum_{j\in\mathbb{Z}} \int_{\mathbb{R}^3} W_0(t,(\xperp, x_{0,n}-2jH),\K) \\
\hspace{2cm}+ W_0(t,(\xperp,-x_{0,n}+2jH),(\Kperp,-k_n)) d\K \\
 \hspace{0.5cm} + \sum_{ j\in \mathbb{Z}^\ast } \int_{\mathbb{R}^3} \e^{- \I k_n \tilde{x}_n} \W_0^{j,\ell=-j} \left(t,  \xperp, x_{0,n}+2jH, \K \right)  d\K, 
\end{multline}
which can be recast as 
\begin{align}
\label{eq:enery_contribution_x_0}
E_{\{x_n = x_{0,n} \}}(t,\xperp,\tilde {x}_n) &= E(t,(\xperp,x_{0,n})) \nonumber \\
&+ 2  \sum_{ j\in \mathbb{Z}^\ast }  \int_{\mathbb{R}^3} \cos ( k_n \tilde{x}_n ) a^{j,\ell=-j}(t, (\xperp,x_{0,n} + 2jH),\K)\DD(\K) d\K \nonumber \\
&+2  \sum_{ j\in \mathbb{Z}^\ast }  \int_{\mathbb{R}^3} \sin ( k_n \tilde{x}_n ) a^{j,\ell=-j}(t, (\xperp,x_{0,n} + 2jH),\K)\widetilde \DD(\K) d\K,
\end{align}
with $\DD$ given by Eq.~\eqref{def:D} and 
\begin{equation}
\label{slab:D_tilde}
\widetilde \DD(\K):=\frac{1}{  c_0 \left|  \K  \right|}
\begin{pmatrix} 
0 & 1\\
-1 & 0
\end{pmatrix}.
\end{equation}  
Finally, using the explicit formulation~\eqref{eq:explicite_form_ajl} for $a^{j\ell}$, the total energy reads 
\begin{align}
\label{eq:enery_contribution_x_0_explicit}
E_{\{x_n = x_{0,n} \}}(t,\xperp,\tilde {x}_n) &= E(t,(\xperp,x_{0,n})) \nonumber \\
&+ 2  \sum_{ j\in \mathbb{Z}^\ast }  \int_{\mathbb{R}^3} \cos ( k_n \tilde{x}_n )  \e^{-\Sigma(\K)t}  \mathbb{A}(\K)\delta(\xperp- \x_{0,\perp}-c_0t\widehat{\K}_\perp) \nonumber \\
&\hspace{2cm} \times \delta(2jH - c_0t\widehat{\K}_n)   \DD(\K) d\K \nonumber \\
&+ 2  \sum_{ j\in \mathbb{Z}^\ast }  \int_{\mathbb{R}^3} \sin ( k_n \tilde{x}_n )  \e^{-\Sigma(\K)t}  \mathbb{A}(\K)  \delta(\xperp- \x_{0,\perp}-c_0t\widehat{\K}_\perp)\nonumber \\
&\hspace{2cm} \times \delta(2jH - c_0t\widehat{\K}_n)  \widetilde \DD(\K) d\K,
\end{align}
where the interference effects are explicitly identified w.r.t. Eq.~\eqref{eq:energy1}. In this formulation, the Dirac mass $\delta(2jH - c_0t\widehat{\K}_n)$ translates that no contribution can be observed before time $2H/c_0$, and their periodic onset occur at each time $2jH/c_0$, for $j\geq 1$. Also, each contributions is damped at rate $\Sigma(\K)$ and freely transport along the plane $\{x_n=x_{0,n}\}$.

Regarding the interference occurring within one wavelength of the plan $\{x_n = H-x_{0,n}\}$, we proceed in a similar way, but with
\begin{multline}
E_{\{x_n = H-x_{0,n}\}}(t,\xperp,\tilde {x}_n) = \sum_{j\in\mathbb{Z}} \int_{\mathbb{R}^3} \lim_{\epsilon\to 0} W_\epsilon(t,(\xperp,  H-x_{0,n}-2jH+ \epsilon \tilde x_n/2),\K)  \\
+ \lim_{\epsilon\to 0}W_\epsilon\left(t,\left(\xperp,-H + x_{0,n}+2jH-\epsilon \tilde x_n/2 \right),(\Kperp,-k_n)\right) d\K \\
+  \frac{1}{2}  \sum_{ j\in \mathbb{Z} } \int_{\mathbb{R}^3} \lim_{\epsilon\to 0} \epsilon  W_\epsilon^{j, \ell = -1-j} (t, (\xperp,-H+x_{0,n}-\epsilon \tilde x_n/2),(\K_\perp,-\epsilon \tilde k_n /2))d\K_\perp d\tilde k_n.
\end{multline}
As to obtain Eq.~\eqref{eq:enery_contribution_x_0_explicit}, the total energy reads this time
\begin{align}
\label{eq:energy_contribution_H-x_0}
E_{\{x_n = H-x_{0,n} \}}(t,\xperp,\tilde {x}_n) &= E(t,(\xperp,H-x_{0,n})) \nonumber \\
&+ 2  \sum_{ j\in \mathbb{Z} }  \int_{\mathbb{R}^3} \cos ( k_n \tilde{x}_n )  \e^{-\Sigma(\K)t} \mathbb{A}(\K)\delta(\xperp-\x_{0,\perp}-c_0t\widehat{\K}_\perp) \nonumber \\
&\hspace{2cm} \times \delta((2j+1)H - c_0t\widehat{\K}_n)   \DD(\K) d\K \nonumber \\
&+ 2  \sum_{ j\in \mathbb{Z} }  \int_{\mathbb{R}^3} \sin ( k_n \tilde{x}_n )  \e^{-\Sigma(\K)t} \mathbb{A}(\K)\delta(\xperp-\x_{0,\perp}-c_0t\widehat{\K}_\perp)\nonumber \\
&\hspace{2cm} \times \delta((2j+1)H - c_0t\widehat{\K}_n)  \widetilde \DD(\K) d\K ,
\end{align}
with an explicit identification of the interferences w.r.t. Eq.~\eqref{eq:energy1}. Again, each contributions is damped at rate $\Sigma(\K)$ and freely transport along the plane $\{x_n=H-x_{0,n}\}$, but the Dirac mass $\delta((2j+1)H - c_0t\widehat{\K}_n)$ tells us that no interference can be observed before time $H/c_0$, and their periodic onsets occur at times $(2j+1)H/c_0$.


\section{Amplifications and weak localization for other sets of boundary conditions} 
\label{sec:other_boundary_conditions}
After considering homogeneous Neumann conditions on both sides of the slab, we consider in this section other sets of boundary conditions, and study the impact on the energy amplifications at the boundaries and within the slab.

\subsection{Dirichlet boundary conditions} 
\label{sec:dirichlet_boundary_condition}

\begin{figure}[tbhp]
\centering \includegraphics[trim = 1.1cm 1.7cm 0.1cm 0.5cm, clip, scale=0.18]{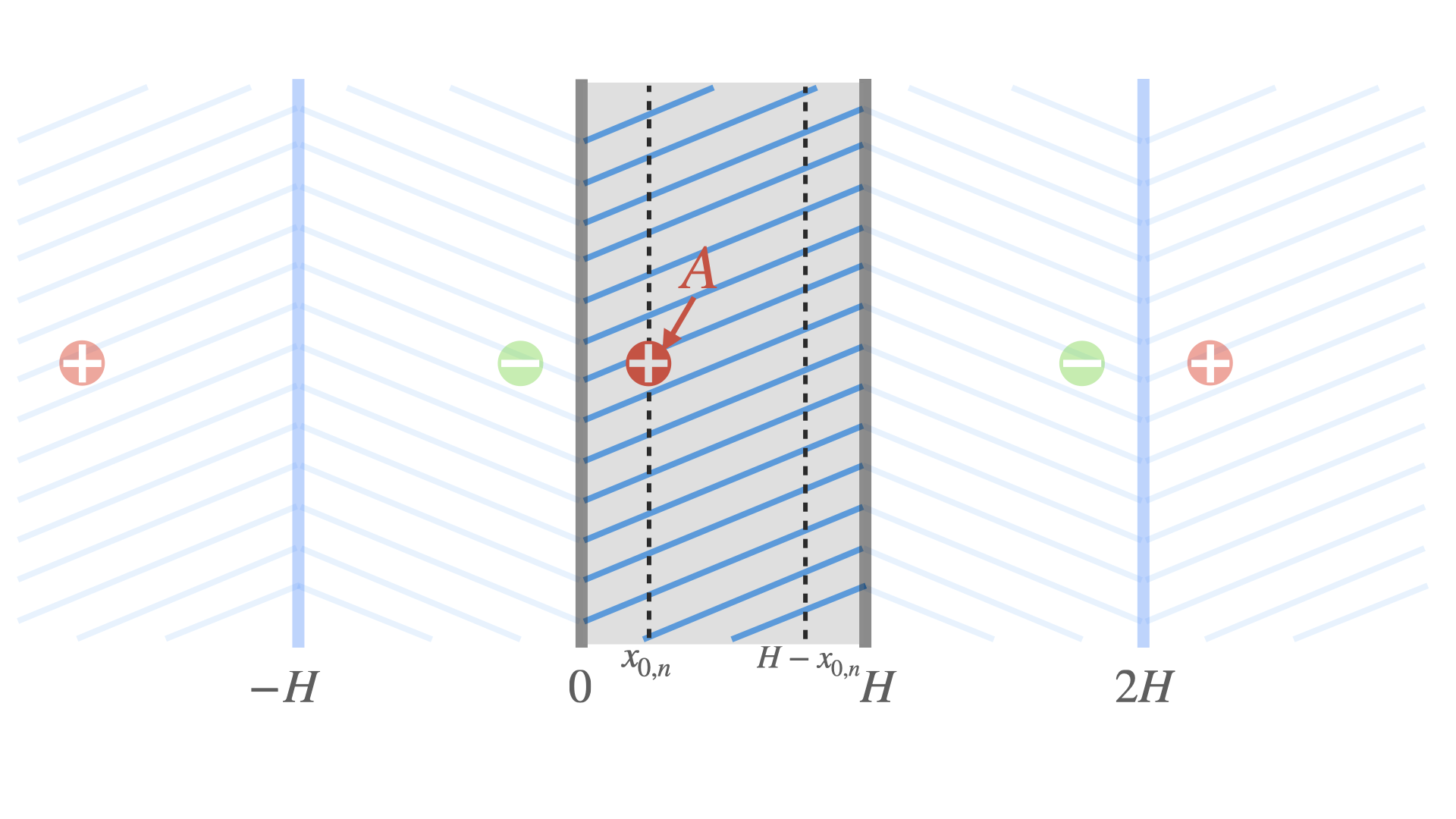}
\caption{\label{fig:superposition_ondes_dirichlet}Sketch of the method of images for the medium delimited by the boundaries $\{\xn=0\}$ and $\{\xn=H\}$, under homogeneous Dirichlet boundary conditions. Each wave propagates in a periodic full-space from an initial condition whose support is centered on one of the red or green symbols. The pluses (in red) indicate a positive initial condition, while the minuses (in green) indicate a negative initial condition (change of sign with respect to the original initial condition). A cancellation of the total energy is observed within one wavelength of the boundaries of the slab (thick grey lines), as well as an amplification of the coherent energy inside the slab, along the dashed lines.}

\end{figure}

In this section, Dirichlet boundary conditions complete the equation~\eqref{eq:equation_des_ondes}. In other words, we consider here 
\begin{equation*}
p(t,\xperp,\xn=0)=p(t,\xperp,\xn=H)=0, \qquad \left(t, \xperp \right) \in \mathbb{R}_+^*\times \mathbb{R}^2.
\end{equation*}
To be compatible with these boundary conditions, the extended wave field $p^\sharp_\epsilon$ is here, instead of  Eq.~\eqref{eq:superposition_H.F}, of the form
\begin{equation}
\label{eq:superposition_H.F_dirichlet_neumann}
p^\sharp_\epsilon(t,\x) =  \sum_{j\in \mathbb{Z}} q_\epsilon(t,\xperp, x_n - 2jH) - q_\epsilon(t,\xperp, -x_n + 2jH), \quad (t,\x) \in  \mathbb{R}_+\times \mathbb{R}^3.
\end{equation}
Each $q_\epsilon(t,\xperp, x_n - 2jH)$ represents a wave propagating in the periodically-structured full-space, whose source location is $\left(\x_{0,\perp},x_{0,n}+2jH\right)$, as labeled with red pluses in Fig.~\ref{fig:superposition_ondes_dirichlet}, and each $-q_\epsilon(t,\xperp, -x_n +2jH )$ represents a wave whose source location is $\left(\x_{0,\perp},-x_{0,n}+2jH\right)$, as labeled with green minuses in Fig.~\ref{fig:superposition_ondes_dirichlet}. The pluses and minuses on the locations of the initial conditions aim to represent the sign with which each wave appears in the decomposition~\eqref{eq:superposition_H.F_dirichlet_neumann} of the total field $p^\sharp_\epsilon$. In the present context, the symmetry relations~\eqref{eq:sym_p-p+} become
\begin{equation*}
p_\epsilon^\sharp(t,\xperp,-\xn) = -p_\epsilon^\sharp(t,\xperp,\xn) \quad \mathrm{and} \quad  p_\epsilon^\sharp(t,\xperp,\xn + 2jH) = p_\epsilon^\sharp(t,\xperp,\xn),       
\end{equation*}
where $(t,\x) \in  \mathbb{R}_+\times \mathbb{R}^3$ and $  j\in\mathbb{Z}$, ensuring that
\begin{equation*}
p_\epsilon^\sharp(t,\xperp,\xn=0)=p_\epsilon^\sharp(t,\xperp,\xn=H)=0, \qquad \left(t, \xperp \right) \in \mathbb{R}_+^*\times \mathbb{R}^2,
\end{equation*}
and
\[
p_\epsilon^\sharp(t,\x) = p_\epsilon(t,\x) \qquad (t,\x)\in \mathbb{R}_+\times\Omega.
\]
This guarantees that $p_\epsilon^\sharp$ is a proper extension to the original wave field $p_\epsilon$.  

In terms of Wigner transform, the decomposition~\eqref{eq:wigner_total} now reads 
\begin{multline}
\label{eq:wigner_total_dirichlet}
W_\epsilon^\sharp(t,\x,\K) = \sum_{j\in \mathbb{Z}} W_\epsilon\left(t,(\xperp,x_n-2jH), \K\right) +  W_\epsilon\left(t,(\xperp,-x_n+2jH),(\Kperp,-k_n) \right)     \\
 -\sum_{j,\ell \in \mathbb{Z}} V_\epsilon^{j\ell}(t,\x,\K) + V_\epsilon^{j\ell}\left(t,(\xperp,-x_n), (\Kperp,-k_n)\right) \\
 + \sum_{\substack{ j,l \in \mathbb{Z}   \\ \ell \ne j}} W_\epsilon^{j\ell}(t,\x,\K) + W_\epsilon^{j\ell}\left(t,(\xperp,-x_n), (\Kperp,-k_n)\right), 
\end{multline}
with a $-$ sign for all the cross-Wigner transforms $V_\epsilon^{j\ell}$ playing a role at the boundaries.

The asymptotic analysis, for each of these Wigner transforms separately, is the same as for Neumann boundaries conditions. Without any focus on specific regions of the slab, the density energy propagation is described by the self-Wigner transforms through Eq.~\eqref{eq:W_0_self}, and despite the change of boundary conditions the reflection conditions~(\ref{eq:reflecting_boundary_conditions0},       
\ref{eq:reflecting_boundary_conditionsH}) are still valid. Regarding the cross-Wigner transforms, the limits~\eqref{eq:asymptotic_limit_cross_wigner_1}, Eq.~\eqref{eq:asymptotic_limit_cross_wigner_1_bis},   Eq.~\eqref{eq:asymptotic_limit_cross_wigner_2}, and Eq.~\eqref{eq:asymptotic_limit_cross_wigner_2_bis} still hold true. Hence, following the same steps as in Sec.~\ref{sec:amplification_at_the_booundaries_for_the_cross_Wigner_transforms}, there is only a change of sign in the total energies within one wavelength of each boundary, which read
\begin{equation}\label{eq:E0_D}
E_{\{x_n = 0\}}(t,\xperp,\tilde {x}_n) = 2 \sum_{j\in \mathbb{Z}}  \int_{\mathbb{R}^3}a\left( t,(\xperp,2jH),\K \right)  (1-\cos(k_n \tilde{x}_n))  \DD(\K) d\K, 
\end{equation}
and 
\begin{equation}\label{eq:EH_D}
E_{\{x_n = H\}}(t,\xperp, \tilde {x}_n) =  2 \sum_{j\in \mathbb{Z}}  \int_{\mathbb{R}^3}  a\left( t,(\xperp,(2j+1)H),\K \right) (1-\cos(k_n \tilde{x}_n))  \DD(\K) d\K,
\end{equation}
where $a$ is the solution to Eq.~\eqref{eq:equation_de_transport}, and $\DD$ is defined by Eq.~\eqref{def:D}. At exactly the boundaries, that is $\tilde {x}_n = 0$ (remembering the change of variable~\eqref{def:chg_var1} for $\ell=-j$ and $\ell = 1-j$ or $-1-j$), we can observe a canceling of the energy intensity at the boundaries
\[E_{\{x_n=0\}}(t,\xperp,\tilde {x}_n = 0) = E_{\{x_n=H\}}(t,\xperp,\tilde {x}_n = 0) = 0,\]
which is consistent with Dirichlet boundary conditions. 

Inside the slab, because of the positive sign in front of all the $W_\epsilon^{j\ell}$ in Eq.~\eqref{eq:wigner_total_dirichlet}, the interference phenomena are exactly the same as the ones described in Sec.~\ref{sec:weak_localization_phenomena_into_the_slab} through Eq.~\eqref{eq:enery_contribution_x_0_explicit} and~\eqref{eq:energy_contribution_H-x_0}.


\subsection{Mixed boundary conditions} 
\label{sec:dirichlet__and_neumann_boundary_conditions}
In this section, the wave equation~\eqref{eq:equation_des_ondes} is completed with a homogeneous Dirichlet boundary condition at $\xn=0$ and a homogeneous Neumann boundary condition at $\xn=H$, that is to say:
\begin{equation*}
p(t,\xperp,\xn=0)=0    \qquad \mathrm{and} \qquad   \partial_{n}p(t,\xperp,\xn=H)=0,
\end{equation*}
that is a Dirichlet condition at $\{x_n=0\}$ and a Neumann condition at $\{x_n=H\}$. The opposite choice can readily be treated with the analysis developed below and straightforward adaptations. In this context, the extended wave field $p^\sharp_\epsilon$ becomes after some algebra
\begin{multline}\label{eq:p_mixed}
p^\sharp_\epsilon(t,\x)  =  \sum_{ j =2k } q_\epsilon(t,\xperp, x_n - 2jH) - \sum_{j =2k+1 } q_\epsilon(t,\xperp, x_n - 2jH) \\
\hspace{0.5cm} - \sum_{j =2k } q_\epsilon(t,\xperp, -x_n + 2jH) + \sum_{j =2k+1 } q_\epsilon(t,\xperp, -x_n + 2jH).
\end{multline}
To obtain this formulation, an even mapping w.r.t $\xn=0$ is first considered, then an odd mapping w.r.t $\xn=H$, and finally an extension by $4H$-periodicity is constructed.  
\begin{figure}[tbhp]
\centering \includegraphics[trim = 1.1cm 1.7cm 0.1cm 0.5cm, clip, scale=0.18]{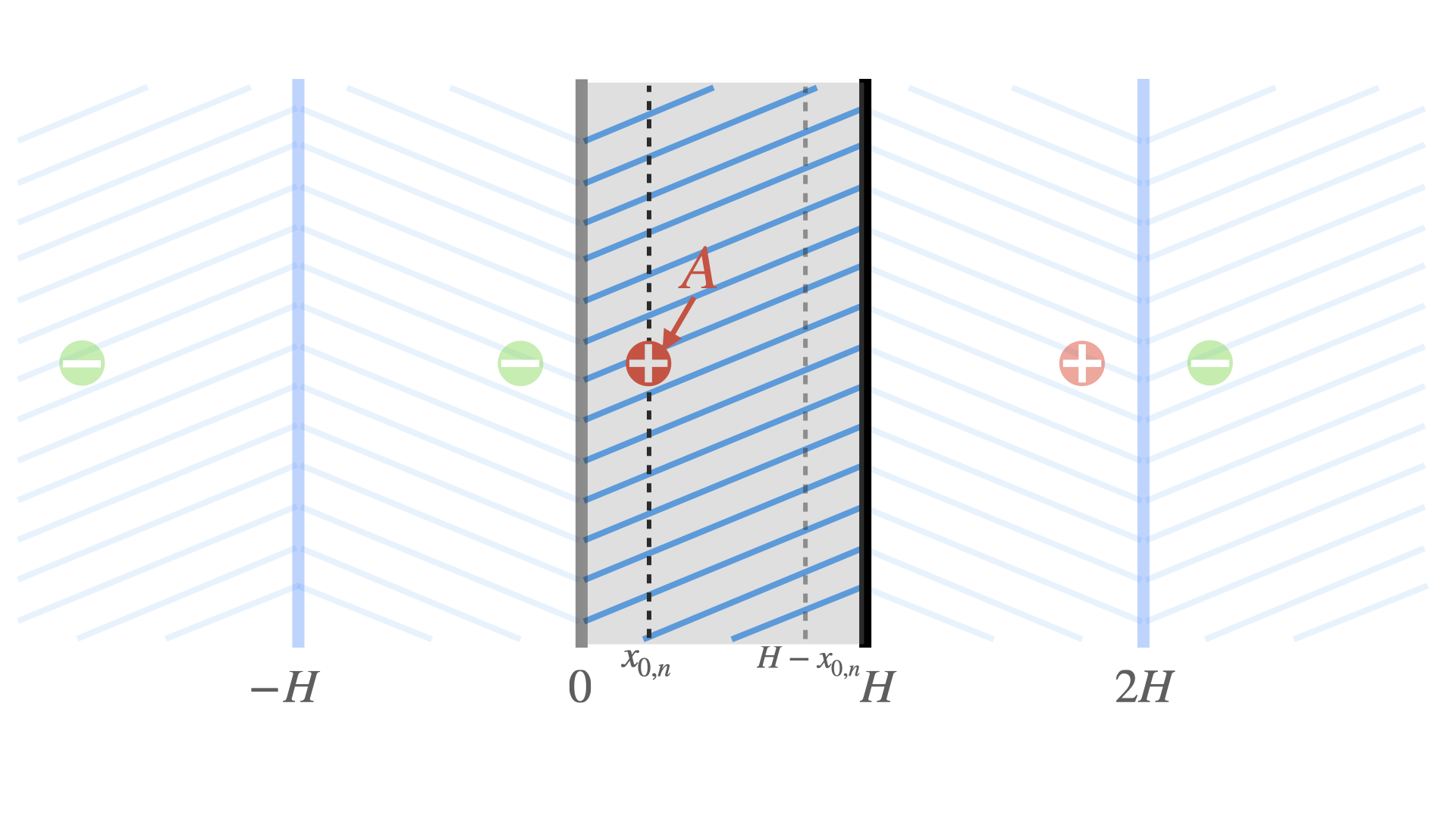}
\caption{\label{fig:superposition_ondes_dirichlet_neumann} Sketch of the method of images for the medium delimited by the boundaries $\{\xn=0\}$ and $\{\xn=H\}$, under mixed boundary conditions. Each wave propagates in a periodic full-space from an initial condition whose support is centered on one of the red or green symbols. The pluses (in red) indicate a positive initial condition, while the minuses (in green) indicate a negative initial condition (change of sign with respect to the original initial condition). A cancellation of the total energy is observed within one wavelength of the Dirichlet boundary (thick grey lines), as well as an amplification of the total energy along the Neumann boundary condition (thick black line), and an amplification of the coherent energy along the dashed lines.}
\end{figure}
As in the previous section, each $q_\epsilon(t,\xperp, x_n - 2jH)$ and $q_\epsilon(t,\xperp, -x_n + 2jH)$ represent waves propagating in the periodically-structured full-space, with supports of their initial conditions centered around $(\x_{0,\perp},x_{0,n}+2jH)$ and $(\x_{0,\perp},-x_{0,n}+2jH)$, respectively, as illustrated in Figure~\ref{fig:superposition_ondes_dirichlet_neumann}. In this picture, the red pluses represent waves whose initial condition corresponds to Eq.~\eqref{eq:cond_init_eps} while the green minuses represent the waves associated to initial conditions with opposite signs.

In the present context, using Eq.~\eqref{eq:p_mixed}, the Wigner transform for the extended wave field now reads
\begin{equation}
\label{eq:wigner_total_dirichlet_neumann}
W_\epsilon^\sharp(t,\x,\K) = W_\epsilon^{1,\sharp}(t,\x,\K)+W_\epsilon^{2,\sharp}(t,\x,\K)+W_\epsilon^{3,\sharp}(t,\x,\K)
\end{equation}
where the three contribution are given by
\begin{multline*}
W_\epsilon^{1,\sharp}(t,\x,\K)  = \sum_{j\in \mathbb{Z}} W_\epsilon\left(t,(\xperp,x_n-2jH), \K\right) \\
+  W_\epsilon\left(t,(\xperp,-x_n+2jH),(\Kperp,-k_n) \right),
\end{multline*}
describing the overall energy density over the whole slab,
\begin{multline}
W_\epsilon^{2,\sharp}(t,\x,\K)  = - \sum_{ j =2k } \sum_{ \ell =2m } V_\epsilon^{j\ell}(t,\x,\K) + V_\epsilon^{j\ell}\left(t,(\xperp,-x_n), (\Kperp,-k_n)\right)  \\
\hspace{0.5cm}+\sum_{ j =2k } \sum_{ \ell =2m+1 } V_\epsilon^{j\ell}(t,\x,\K) + V_\epsilon^{j\ell}\left(t,(\xperp,-x_n), (\Kperp,-k_n)\right)  \\
\hspace{0.5cm}+\sum_{ j =2k+1 } \sum_{ \ell =2m } V_\epsilon^{j\ell}(t,\x,\K) + V_\epsilon^{j\ell}\left(t,(\xperp,-x_n), (\Kperp,-k_n)\right)  \\
\hspace{0.5cm}-\sum_{ j =2k+1 } \sum_{ \ell =2m+1 } V_\epsilon^{j\ell}(t,\x,\K) + V_\epsilon^{j\ell}\left(t,(\xperp,-x_n), (\Kperp,-k_n)\right),
\end{multline}
allowing to describe the interferences taking place at the boundaries $\{x_n=0\}$ and $\{x_n=H\}$, and  
\begin{multline}
W_\epsilon^{3,\sharp}(t,\x,\K) = \sum_{ j =2k }\sum_{ \substack{ \ell =2m \\ m \ne k}} W_\epsilon^{j\ell}(t,\x,\K) + W_\epsilon^{j\ell}\left(t,(\xperp,-x_n), (\Kperp,-k_n)\right)  \\
\hspace{0.5cm}-\sum_{ j =2k } \sum_{ \ell =2m+1 }  W_\epsilon^{j\ell}(t,\x,\K) + W_\epsilon^{j\ell}\left(t,(\xperp,-x_n), (\Kperp,-k_n)\right) \\
\hspace{0.5cm}-\sum_{ j =2k+1 } \sum_{ \ell =2m }  W_\epsilon^{j\ell}(t,\x,\K) + W_\epsilon^{j\ell}\left(t,(\xperp,-x_n), (\Kperp,-k_n)\right) \\
\hspace{0.5cm}+\sum_{ j =2k+1 } \sum_{\substack{ \ell =2m+1 \\ m \ne k} }  W_\epsilon^{j\ell}(t,\x,\K) + W_\epsilon^{j\ell}\left(t,(\xperp,-x_n), (\Kperp,-k_n)\right),
\end{multline}
allowing to describe the interferences taking place within the slab at $\{x_n=x_{0,n}\}$ and $\{x_n=H-x_{0,n}\}$.

As for the previous cases, the contribution to $W_\epsilon^{1,\sharp}$ is given by Eq.~\eqref{eq:limW1} to describe the energy propagation over the whole slab. Regarding the contribution to $W_\epsilon^{2,\sharp}$ at the boundaries, following the same steps as in Sec.~\ref{sec:contribution_v_eps} and Sec.~\ref{sec:amplification_at_the_booundaries_for_the_cross_Wigner_transforms} yields necessarily to the choices $\ell=-j$ (both indexes have the same parity) at $\{x_n = 0\}$ and $\ell=1-j$ or $\ell=-1-j$ (one even, one odd) at $\{x_n = H\}$. As a result, within one wavelength of the boundaries the wave energy reads
\begin{equation}\label{eq:E0_M}
E_{\{x_n = 0\}}(t,\xperp,\tilde {x}_n) = 2 \sum_{j\in \mathbb{Z}}  \int_{\mathbb{R}^3}a\left( t,(\xperp,2jH),\K \right)  (1-\cos(k_n \tilde{x}_n))  \DD(\K) d\K, 
\end{equation}
and
\begin{equation}\label{eq:EH_M}
E_{\{x_n = H\}}(t,\xperp, \tilde {x}_n) =   2 \sum_{j\in \mathbb{Z}}  \int_{\mathbb{R}^3}  a\left( t,(\xperp,(2j+1)H),\K \right) (1+\cos(k_n \tilde{x}_n))  \DD(\K) d\K. 
\end{equation}
If we look exactly at the boundary, that is $\tilde x_n =0$, we observe again a canceling of the energy intensity at $\{x_n = 0\}$ (with a Dirichlet condition)
\[E_{\{x_n=0\}}(t,\xperp,\tilde {x}_n = 0) = 0,\]
and a doubling of the energy w.r.t. Eq.~\eqref{eq:energy1} at $\{x_n = H\}$ (with a Neumann condition)
\[E_{\{x_n = H\}}(t,\xperp,\tilde {x}_n = 0) = 2 \, E(t,(\xperp,H)).\]

Within the slab, the interference effects are provided by $W_\epsilon^{3,\sharp}$. Following the steps of Sec.~\ref{sec:contribution_W_eps} and Sec.~\ref{sec:weak_localization_phenomena_into_the_slab}, we necessarily have $\ell =-j$ (both indexes have the same parity) at $\{x_n = x_{0,n}\}$ and $\ell=-1-j$ (one even, one odd) at $\{x_n = H-x_{0,n}\}$. Hence, at the plane $ \{x_n = x_{0,n} \}$, the weak localization phenomena described in Sec.~\ref{sec:weak_localization_phenomena_into_the_slab} and Eq.~\eqref{eq:enery_contribution_x_0_explicit} still holds true with an increase of the energy w.r.t. Eq.~\eqref{eq:energy1}. However, in the vicinity of the plane $ \{x_n = H-x_{0,n} \} $, Eq.~\eqref{eq:energy_contribution_H-x_0} becomes now 
\begin{align}
E_{\{x_n = H-x_{0,n} \}}(t,\xperp,\tilde {x}_n) &= E(t,(\xperp,H-x_{0,n})) \nonumber\\
& - 2  \sum_{ j\in \mathbb{Z} }  \int_{\mathbb{R}^3} \cos ( k_n \tilde{x}_n ) \e^{-\Sigma(\K)t} \mathbb{A}(\K)\delta\big(\xperp-\x_{0\perp}-c_0t\widehat{\K}_\perp\big) \nonumber \\
& \hspace{2cm}\delta\big((2j+1)H - c_0t\widehat{k}_n\big)   \DD(\K) d\K \nonumber \\
& - 2  \sum_{ j\in \mathbb{Z} }  \int_{\mathbb{R}^3} \sin( k_n \tilde{x}_n ) \e^{-\Sigma(\K)t} \mathbb{A}(\K)\delta\big(\xperp-\x_{0\perp}-c_0t\widehat{\K}_\perp\big) \nonumber \\
& \hspace{2cm} \delta\big((2j+1)H - c_0t\widehat{k}_n\big)\widetilde \DD(\K) d\K,
\end{align}
yielding a decrease of the total energy at exactly $x_n=H-x_{0,n}$ for $\tilde x_n=0$.


\section*{Conclusion}

Radiative transfer equations have been derived for acoustic waves propagating in a randomly-fluctuating slab in the weak-scattering regime. Our approach has been built on an asymptotic analysis of the Wigner transform of the wave solution and the method of images. More specifically, the method of images consists in representing the solution of the wave equation~\eqref{eq:equation_des_ondes}, supported on the slab, with a superposition of several solutions of the same wave equation~\eqref{eq:equation_des_onde_full_space} extended to a full-space with a periodic map of mechanical properties. In this superposition, one term is equipped with the original initial condition, and the others with initial conditions chosen so as to enforce the boundary conditions. The energy densities for these wave fields verify the same RTE~\eqref{eq:equation_de_transport} as in an unbounded domain, with reflection conditions~\eqref{eq:reflecting_boundary_conditions} at the boundaries that are the same whether the original wave equation was completed with Dirichlet or Neumann boundary conditions. Additionally, due to interference effects, under Neumann boundaries conditions, amplifications of the total energy in the vicinity of the boundaries and the planes $\{x_n=x_{0,n}\}$ and $\{x_n=H-x_{0,n}\}$ are observed. Under Dirichlet boundary conditions, the same amplifications inside the domain are observed, while cancellation of the energy take place near the boundaries. 
When the interface $\{x_n=0\}$ is equipped with a Dirichlet boundary condition and the other interface $\{x_n=H\}$ with a Neumann boundary condition, interference effects lead to amplifications of the total energy in the vicinity of the planes $\{ x_n=H\}$ and $\{x_n=x_{0,n}\}$, while decreases in the vicinity of the planes $\{x_n=0\}$ and $\{x_n=H-x_{0,n}\}$ are observed. The amplifications along the boundaries apply to the entire energy density field while the modulation inside the domain only applies to the coherent field.

\appendix

\section{Spectral decomposition of the dispersion matrix}
\label{sec:appendix_spectral_decomposition}

This section is devoted to the derivation of the decomposition~\eqref{eq:W_0_self} following the same steps and uses the same formalism as described in~\cite{Bal2005} for two signals propagating in two different media. To this end, let us consider the Wigner transform $\W_\epsilon = W[\mathbf{g}_\epsilon(t,\cdot),\p_\epsilon(t,\cdot)]$, where $\mathbf{g}_\epsilon$ and $\p_\epsilon$ are two vector fields satisfying the $2\times2$ first order hyperbolic systems :
\begin{equation}
\label{eq:equation_hyperbolique_ordre_1_1}
\epsilon\partial_t \mathbf{g}_\epsilon + \mathcal{A}_\epsilon^1\mathbf{g}_\epsilon= 0, \qquad \mathrm{ where }\qquad \mathcal{A}_\epsilon^1:=-\begin{pmatrix}0& c_{1}^2 \left(\frac{\x}{\epsilon}\right)\\ \epsilon^2\Delta & 0  \end{pmatrix}, 
\end{equation}
and 
\begin{equation}
\label{eq:equation_hyperbolique_ordre_1_2}
\epsilon\partial_t \p_\epsilon+ \mathcal{A}_\epsilon^2\p_\epsilon = 0, \qquad \mathrm{ where }\qquad \mathcal{A}_\epsilon^2:=-\begin{pmatrix}0& c_{2}^2 \left(\frac{\x}{\epsilon}\right)\\ \epsilon^2\Delta & 0  \end{pmatrix}. 
\end{equation}
According to Eq.~\eqref{eq:sound_speed}, the operator $\mathcal{A}_\epsilon^\phi$ can be split as
\begin{equation*}
\mathcal{A}_\epsilon^\phi=-\begin{pmatrix}0 & c_0^2 \\  \epsilon^2 \Delta& 0  \end{pmatrix} + \sqrt{\epsilon}\, \nu_\phi \left(\x, \frac{\x}{\epsilon} \right)K, \qquad \text{where}\qquad K =\begin{pmatrix}0& 1\\ 0& 0  \end{pmatrix},\quad \phi = 1,2.
\end{equation*}
exhibiting the separation between the slow and fast components of the medium fluctuations introduced in Sec.~\ref{sec:structure_of_the_inhomogeneities}. From Eq.~\eqref{eq:equation_hyperbolique_ordre_1_1} and Eq.~\eqref{eq:equation_hyperbolique_ordre_1_2} we deduce the following relation for $\W_\epsilon^{j\ell}$
\begin{equation}
\label{eq:eq_wigner_annex}
\epsilon \partial_t \W_\epsilon(t,\x,\K)  + W\left[\mathcal{A}_\epsilon^1 \mathbf{g}_\epsilon,\p_\epsilon \right](t,\x,\K) \\
+ W\left[\mathbf{g}_\epsilon,\mathcal{A}_\epsilon^2 \p_\epsilon \right](t,\x,\K) =0,
\end{equation}
In what follows, the Laplacian operator, $\epsilon^2 \Delta $ in $\mathcal{A}_\epsilon^\phi$, is rewrite as the pseudo-differential operator $\F\left( \epsilon \mathbf{D} \right)$ defined by
\begin{equation}
\label{def:pseudo_diff}
\F\left( \epsilon \mathbf{D} \right)\left[\U\right](\x)=\int_{\mathbb{R}^3} \e^{\I\p \cdot \x} \left(\I\epsilon\p \right) \cdot \left(\I\epsilon\p \right) \widehat{\U}(\p)\frac{d\p}{(2\pi)^3}.
\end{equation}
Knowing the four following relations  
\begin{align*}
W\left[\F(\epsilon \mathbf{D})\mathbf{g}_\epsilon,\p_\epsilon \right](t,\x,\K) & = \F\left( \I\K + \frac{\epsilon \mathbf{D}}{2} \right)\left[\W_\epsilon\right](t,\x,\K),\\
W\left[\mathbf{g}_\epsilon,\F(\epsilon \mathbf{D})\p_\epsilon\right](t,\x,\K) & = \F\left( \I\K - \frac{\epsilon \mathbf{D}}{2} \right)\left[\W_\epsilon\right](t,\x,\K),\\
W\left[\nu_1\left( \x,\frac{\x}{\epsilon} \right)\mathbf{g}_\epsilon,\p_\epsilon\right](t,\x,\K) & = \int_{\mathbb{R}^3}\e^{\I \x \cdot\p/\epsilon}\widehat{\nu}_1(\x,\p)\W_\epsilon\left(t,\x,\K-\frac{\p}{2}\right)\frac{d\p}{(2\pi)^3} \\
& + O\left( \epsilon \right),\\
W\left[\mathbf{g}_\epsilon,\nu_2\left( \x,\frac{\x}{\epsilon} \right)\p_\epsilon\right](\x,\K)& = \int_{\mathbb{R}^3}\e^{\I \x \cdot\p/\epsilon}\widehat{ \nu}_2 (\x,\p)\W_\epsilon\left(t,\x,\K+\frac{\p}{2}\right)\frac{d\p}{(2\pi)^3} \\
& + O\left( \epsilon \right) ,
\end{align*}
the relation~\eqref{eq:eq_wigner_annex} can be recast as
\begin{multline}
\label{eq:eq_wigner_annex_recast}
\epsilon \partial_t \W_\epsilon (t,\x,\K)  + \Pdiff\left(\I\K + \frac{\epsilon\mathbf{D}_\x}{2}    \right)\W_\epsilon (t,\x,\K) 
 + \W_\epsilon(t,\x,\K) \Pdiff^\ast \left(\I\K  - \frac{\epsilon\mathbf{D}_\x}{2}    \right)\\
+ \sqrt{\epsilon}\left( K \, \mathcal{K}_\epsilon^{1-}  [\W_\epsilon](t,\x,\K) +  \mathcal{K}^{2+}_\epsilon[\W_\epsilon](t,\x,\K) \, K^\ast \right) =0,
\end{multline}
where 
\begin{equation}
\mathcal{K}^{\phi \pm}_\epsilon [W](t,\x,\K)=\int_{\mathbb{R}^3}e^{i\x\cdot\p/\epsilon}\widehat{\nu}_\phi(\x,\p)W\left(t,\x,\K\pm\frac{\p}{2}\right)\frac{d\p}{(2\pi)^3}, \qquad \phi = 1,2,
\end{equation}
and
\begin{equation}\label{def:P}
\mathcal{P}\left(\I \K + \frac{\epsilon \mathbf{D}}{2} \right) := -\begin{pmatrix}
0 & c_0^2 \\ \mathcal{F}\left(\I \K + \frac{\epsilon \mathbf{D}}{2} \right) & 0 
\end{pmatrix}.
\end{equation}
Due to the presence of the rapidly oscillating phases $\e^{\I \x \cdot\p/\epsilon}$ in $\mathcal{K}^{\phi\pm}_\epsilon$, the fast variable $\y=\x /\epsilon$ is introduced, and $\W_\epsilon$ is rewritten as
\begin{equation*}
\W_\epsilon(t,\x,\K)  = \W_\epsilon(t,\x,\y,\K)_{|\y=\x /\epsilon}=\W_\epsilon\left(t,\x,\frac{\x}{\epsilon},\K\right),
\end{equation*}   
to account for this new variable. Having two spatial variables for $\W_\epsilon$, the differential operator $\mathbf{D}$ is now given by
\[
\mathbf{D} = \mathbf{D}_\x + \frac{1}{\epsilon}\mathbf{D}_\y, 
\]
and Eq.~\eqref{eq:eq_wigner_annex_recast} can be rewritten as
\begin{multline}
\label{eq:wigner_equation_rescale_annex}
\epsilon \partial_t \W_\epsilon (t,\x,\y,\K)  + \Pdiff\left(\I\K+ \frac{\mathbf{D}_\y}{2} + \frac{\epsilon\mathbf{D}_\x}{2}    \right)\W_\epsilon (t,\x,\y,\K) \\
 + \W_\epsilon (t,\x,\y,\K) \Pdiff^\ast \left(\I\K - \frac{\mathbf{D}_\y}{2} - \frac{\epsilon\mathbf{D}_\x}{2}    \right)\\
+ \sqrt{\epsilon}\left( K \, \mathcal{K}^{1-}  [\W_\epsilon](t,\x,\y,\K) +  \mathcal{K}^{2+} [\W_\epsilon](t,\x,\y,\K) \, K^\ast \right) +O(\epsilon^{3/2})=0,
\end{multline}
with
\begin{equation*}
\mathcal{K}^{\phi \pm} [W](t,\x,\y,\K):=\int_{\mathbb{R}^3}e^{i\y\cdot\p}\widehat{\nu}_\phi(\x,\p)W\left(t,\x,\y,\K \pm \frac{\p}{2}\right)\frac{d\p}{(2\pi)^3}, \qquad \phi =1,2.
\end{equation*}
To extract the decomposition~\eqref{eq:W_0_self} from Eq.~\eqref{eq:wigner_equation_rescale_annex}, we consider the following expansion for $\W_\epsilon$ in powers of $\epsilon$
\begin{equation}
\label{eq:asymptotic_expansion_1_annex}
\W_\epsilon(t,\x,\y,\K) = \W_0(t,\x,\K) + \sqrt{\epsilon}\, \W_1(t,\x,\y,\K) + \epsilon \, \W_2(t,\x,\y,\K),
\end{equation} 
so that the asymptotic behavior of $\W_\epsilon$ is characterized by $\W_0$. We also consider the first order expansion
\begin{equation}
\label{eq:asymptotic_expansion_0_annex}
\F\left(\I\K+ \frac{\mathbf{D}_\y}{2} + \frac{\epsilon\mathbf{D}_\x}{2}    \right) = \F\left(\I\K+ \frac{\mathbf{D}_\y}{2}  \right) + \frac{\epsilon}{2}\, \F'\left(\I\K+ \frac{\mathbf{D}_\y}{2}  \right) \cdot\nabla_\x + O(\epsilon^2),
\end{equation}
where the symbol of $\F'$ is given by
\begin{equation*}
\F'(\I \K) := -2\K = 2\I q_0(\I \K)\nabla_\K q_0(\I\K),
\end{equation*}
and the one of $q_0$ is defined as  
\begin{equation*}
q_0(\I \K):=\sqrt{-\mathcal{F}(\I \K)} = |\K|.
\end{equation*}
According to Eq.~\eqref{def:P}, $\Pdiff$ also admits an expansion of the form 
\begin{equation}
\label{eq:asymptotic_expansion_annex}
\Pdiff=\Pdiff_0 + \epsilon\, \Pdiff_1 + \mathrm{O}(\epsilon^2),
\end{equation} 
where each term is defined according to Eq.~\eqref{eq:asymptotic_expansion_0_annex}. Injecting Eq.~\eqref{eq:asymptotic_expansion_1_annex} and Eq.~\eqref{eq:asymptotic_expansion_annex} into Eq.~\eqref{eq:wigner_equation_rescale_annex} yields a sequence of three equations by equating the coefficients associated to each power of~$\epsilon$.

The leading order terms yield the relation
\begin{equation}
\label{def:L0}
\mathcal{L}_0 \W_0 := \Pdiff_0(\I \K) \W_0 + \W_0 \Pdiff_0^*(\I \K)=0.
\end{equation}
In this equation, the dispersion matrix is defined by 
\[
\Pdiff_0(\I \K) = - \begin{pmatrix}
0 & c_0^2 \\
\F(\I \K) & 0
\end{pmatrix}
\]
and admits the following spectral representation
\begin{equation*}
\Pdiff_0 = \lambda_+\B_+\C_+^\ast + \lambda_-\B_-\C_-^\ast, 
\end{equation*}
where
\begin{multline}
\label{def:c+}
\lambda_\pm(\K) := \pm  \I c_0 q_0(\I\K), \qquad \B_\pm(\K):= \frac{1}{\sqrt{2}}\begin{pmatrix} \pm \I q_0^{-1}(\I \K) \\ c_0^{-1} \end{pmatrix} \\ \mathrm{and} \quad \C_\pm(\K):= \frac{1}{\sqrt{2}}\begin{pmatrix} \pm \I q_0(\I \K) \\ c_0 \end{pmatrix},
\end{multline} 
with 
\[\B_\pm^\ast \C_\pm=1.\]
Using that $\left(\B_+(\K), \B_-(\K)\right)$ forms a basis of $\mathbb{R}^2$, the matrix $\W_0$ itself can be decomposed as  
\begin{equation}
\label{def:a+}
\W_0 = \sum_{i,m=\pm}\ed{im} \B_i\B_m^* \qquad\text{with}\qquad \ed{im} := \C_i^\ast \W_0 \C_m. 
\end{equation}
Plugging this relation into Eq.~\eqref{def:L0} gives $\ed{im}=0$ for $i\neq m$, and then
\begin{equation*}
\W_0=\ed+\B_+\B_+^* +\ed- \B_-\B_-^*,\qquad\text{with}\qquad \ed+:= \ed{++}\qquad\text{and}\qquad \ed -:= \ed{- -}.
\end{equation*}
Also, using that $\C (\K)= \C (-\K)$ for any $\K$, we have
\begin{equation}\label{eq:rel_sym_a}
\ed\pm(\K)= \ed\mp(-\K),
\end{equation}
so that we just have to focus our attention on $\ed +$. Keeping in mind the relation~\eqref{eq:rel_sym_a} we drop the $+$ dependence for $a$ and we can deduce when $c_1=c_2=\ctilde$
\begin{equation*}
\W_0(t,\x,\K) := \lim_{\epsilon \to 0} W_\epsilon(t,\x,\K)
= a (t,\x,\K)\BB(\K) + a (t,\x,-\K)\BB^T(\K),
\end{equation*}
where 
\begin{equation*}
\BB(\K):= \frac{1}{2}\begin{pmatrix}
1/\vert \K\vert^2 & i/(c_0\vert \K\vert) \\
- i/(c_0\vert \K\vert) & 1/c_0^2
\end{pmatrix}.
\end{equation*}
As described in \cite{Messaoudi2022}, this analysis can be pushed further by analyzing the role played by the correctors $\W_1$ and $\W_2$ in Eq. \eqref{eq:asymptotic_expansion_1_annex}, and show that $a$ satisfies the RTE~\eqref{eq:equation_de_transport} with the initial condition~\eqref{eq:a0}.

\section{Derivation of Eq.~\eqref{eq:RTE_coherent} }

\label{sec:appendix}

This section is devoted to the derivation of Eq.~\eqref{eq:RTE_coherent} through an asymptotic analysis of the Wigner transform $\W_\epsilon^{j\ell} = W[\mathbf{g}_\epsilon^{j\ell}(t,\cdot),\p_\epsilon^{j\ell}(t,\cdot)]$, where $\mathbf{g}_\epsilon^{j\ell}$ and $\p_\epsilon^{j\ell}$ are given by Eq.~\eqref{eq:cross_waves_vectors_axis}. Knowing that $\mathbf{g}_\epsilon^{j\ell}$ and $\p_\epsilon^{j\ell}$ are two signals propagating in two different media, this analysis follows the same steps as described in~\cite{Bal2005} and some of the results presented in the previous Appendix \ref{sec:appendix_spectral_decomposition} will be also used. In this way,  Eq.~\eqref{eq:equation_hyperbolique_ordre_1_1} and Eq.~\eqref{eq:equation_hyperbolique_ordre_1_2} becomes
\begin{equation}
\epsilon\partial_t \mathbf{g}_\epsilon^{j\ell} + \mathcal{A}_\epsilon^1\mathbf{g}_\epsilon^{j\ell} = 0, \qquad \mathrm{ where }\qquad \mathcal{A}_\epsilon^1:=-\begin{pmatrix}0& \ctilde^2 \left(\frac{\xperp}{\epsilon},\frac{\xn - 2x_{0,n}  - (j+\ell)H}{\epsilon} \right)\\ \F\left( \epsilon \mathbf{D} \right) & 0  \end{pmatrix}, 
\end{equation}
and 
\begin{equation}
\epsilon\partial_t \p_\epsilon^{j\ell}+ \mathcal{A}_\epsilon^2\p_\epsilon^{j\ell} = 0, \qquad \mathrm{ where }\qquad \mathcal{A}_\epsilon^2:=-\begin{pmatrix}0& \ctilde^2 \left(\frac{\xperp}{\epsilon},\frac{\xn - (j+\ell)H}{\epsilon} \right)\\ \F\left( \epsilon \mathbf{D} \right) & 0  \end{pmatrix}. 
\end{equation}
Here $\F\left( \epsilon \mathbf{D} \right)$ is defined by Eq.~\eqref{def:pseudo_diff} and the operator $\mathcal{A}_\epsilon^\phi$ will be split as
\begin{equation*}
\mathcal{A}_\epsilon^\phi=-\begin{pmatrix}0 & c_0^2 \\  \F\left( \epsilon \mathbf{D} \right)& 0  \end{pmatrix} + \sqrt{\epsilon}\, \nu_\phi \left(\x, \frac{\x}{\epsilon} \right)K, \qquad \text{where}\qquad K =\begin{pmatrix}0& 1\\ 0& 0  \end{pmatrix},\quad \phi = 1,2.
\end{equation*}
Here, we have
\begin{align}
\label{def:random_fluctuations_annex}
\nu_1\left(\x,\frac{\x}{\epsilon}\right) & = \nu_\sharp \left((\xperp,\xn - 2x_{0,n}  - (j+\ell)H), \left(\frac{\xperp}{\epsilon}, \frac{\xn - 2x_{0,n}  - (j+\ell)H}{\epsilon}  \right) \right), \nonumber\\
&\\
 \nu_2\left(\x,\frac{\x}{\epsilon}\right)& = \nu_\sharp \left( (\xperp,\xn - (j+\ell)H),\left(\frac{\xperp}{\epsilon}, \frac{\xn  - (j+\ell)H}{\epsilon}  \right)  \right),\nonumber
\end{align} 
where $ \nu_\sharp$ (given by Eq.~\eqref{def:N}) corresponds to the random fluctuations associated to the extended velocity field $c_\sharp$. 
The relation~\eqref{eq:wigner_equation_rescale_annex} can be recast as
\begin{multline}
\label{eq:wigner_equation_rescale_annex_B}
\epsilon \partial_t \W_\epsilon^{j\ell} (t,\x,\y,\K)  
+ \Pdiff\left(\I\K+ \frac{\mathbf{D}_\y}{2} + \frac{\epsilon\mathbf{D}_\x}{2}    \right)\W_\epsilon^{j\ell} (t,\x,\y,\K) \\
+ \W_\epsilon^{j\ell} (t,\x,\y,\K) \Pdiff^\ast \left(\I\K - \frac{\mathbf{D}_\y}{2} - \frac{\epsilon\mathbf{D}_\x}{2}    \right)\\
+ \sqrt{\epsilon}\left( K \, \mathcal{K}^{1-}  [\W_\epsilon^{j\ell}](t,\x,\y,\K) +  \mathcal{K}^{2+} [\W_\epsilon^{j\ell}](t,\x,\y,\K) \, K^\ast \right) +O(\epsilon^{3/2})=0,
\end{multline}
with
\begin{equation*}
\mathcal{K}^{\phi \pm} [W](t,\x,\y,\K):=\int_{\mathbb{R}^3}e^{i\y\cdot\p}\widehat{\nu}_\phi (\x,\p)W\left(t,\x,\y,\K \pm \frac{\p}{2}\right)\frac{d\p}{(2\pi)^3}, \qquad \phi =1,2
\end{equation*}
and $\mathcal{P}$ is given by Eq.~\eqref{def:P}.
To extract the transport equation from Eq.~\eqref{eq:wigner_equation_rescale_annex_B}, we consider the following expansion for $\W_\epsilon^{j\ell}$ in powers of $\epsilon$
\begin{equation}
\label{eq:asymptotic_expansion_1_annex_B}
\W_\epsilon^{j\ell}(t,\x,\y,\K) = \W_0^{j\ell}(t,\x,\K) + \sqrt{\epsilon}\, \W_1^{j\ell}(t,\x,\y,\K) + \epsilon \, \W_2^{j\ell}(t,\x,\y,\K),
\end{equation} 
so that the asymptotic behavior of $\W_\epsilon^{j\ell}$ is characterized by $\W_0^{j\ell}$. 
$\Pdiff$ also admits an expansion of the form 
\begin{equation}
\label{eq:asymptotic_expansion_annex_B}
\Pdiff=\Pdiff_0 + \epsilon\, \Pdiff_1 + \mathrm{O}(\epsilon^2),
\end{equation} 
where each term is defined according to Eq.~\eqref{eq:asymptotic_expansion_0_annex}. Injecting Eq.~\eqref{eq:asymptotic_expansion_1_annex_B} and Eq.~\eqref{eq:asymptotic_expansion_annex_B} into Eq.~\eqref{eq:wigner_equation_rescale_annex_B} yields a sequence of three equations by equating the coefficients associated to each power of~$\epsilon$.

\subsection{Leading order}
\label{sec:leading_order}
The leading order terms yield the relation
\begin{equation*}
\mathcal{L}_0 \W_0^{j\ell} := \Pdiff_0(\I \K) \W_0^{j\ell} + \W_0^{j\ell} \Pdiff_0^*(\I \K)=0,
\end{equation*}
with dispersion matrix $\Pdiff_0(\I \K)$ defined as
\[
\Pdiff_0(\I \K) = - \begin{pmatrix}
0 & c_0^2 \\
\F(\I \K) & 0
\end{pmatrix}.
\]
As a result, the spectral analysis provided in the previous Appendix \ref{sec:appendix_spectral_decomposition} leads to the following decomposition of $\W_0^{j\ell}$: 
\begin{multline}
\label{eq:dec_W0}
\W_0^{j\ell} =\ed+^{j\ell} \B_+\B_+^* +\ed-^{j\ell} \B_-\B_-^*,\qquad\text{with}\qquad \ed+^{j\ell}:=  \C_+^\ast \W_0^{j\ell} \C_+,\\ 
\text{and}\qquad \ed -^{j\ell}:= \C_-^\ast \W_0^{j\ell} \C_-.
\end{multline}
Here, the $\B_\pm$ and $\C_\pm$ are given by Eq.~\eqref{def:c+}. Also, using that $\C_\pm (\K)= \C_\pm (-\K)$ for any $\K$, we have
\begin{equation}\label{eq:rel_sym_ajl}
\ed\pm^{j\ell}(\K)= \ed\mp^{j\ell}(-\K),
\end{equation}
so that we just have to focus our attention on $\ed +^{j\ell}$. Keeping in mind the relation~\eqref{eq:rel_sym_ajl}, we drop the $+$ dependence for $a^{j\ell}$ in the remaining of the proof.

\subsection{First order term $\W_1^{j\ell}$}
\label{sec:first_order}

Equating like powers of $\epsilon^{1/2}$ in Eq.~\eqref{eq:wigner_equation_rescale_annex_B}, together with Eq.~\eqref{eq:asymptotic_expansion_1_annex_B} and Eq.~\eqref{eq:asymptotic_expansion_0_annex}, we obtain the following relation
\begin{equation*}
\Pdiff_0\left(\I\K + \frac{\mathbf{D}_\y}{2}  \right )\W_1^{j\ell} + \W_1^{j\ell} \Pdiff_0^*\left(\I\K - \frac{\mathbf{D}_\y}{2}  \right ) +  \mathcal{K}^{1-}  K \W_0^{j\ell} +  \mathcal{K}^{2+} \W_0^{j\ell} K^\ast = 0.
\end{equation*}
To avoid singular terms and to preserve causality, a regularization term $\theta$ is added following \cite{Bal2005, Ryzhik1996} and will be sent to $0$ later on: 
\begin{equation}
\label{eq:first_order_corrector}
\Pdiff_0\left(\I\K + \frac{\mathbf{D}_\y}{2}  \right )\W_1^{j\ell} + \W_1^{j\ell} \Pdiff_0^*\left(\I\K - \frac{\mathbf{D}_\y}{2}  \right )  + \theta \W_1^{j\ell}+  \mathcal{K}^{1-} K \W_0^{j\ell} + \mathcal{K}^{2+} \W_0^{j\ell} K^\ast  = 0.
\end{equation}
Now, taking the Fourier transform of Eq.~\eqref{eq:first_order_corrector} in $\y$ leads to
\begin{multline}
\label{eq:fourier_first_order_corrector}
\Pdiff_0\left(\I\K + \I\frac{\p}{2}  \right )\widehat{\W}_1^{j\ell} +  \widehat{\W}_1^{j\ell}\Pdiff_0^*\left(\I \K - \I \frac{\p}{2}  \right ) + \theta \widehat{\W}_1^{j\ell}  \\
+  \widehat{\nu}_1 (\x,\p) K \W_0^{j\ell}\left( \K + \frac{\p}{2}   \right) \\
+ \widehat{\nu}_2(\x,\p) \W_0^{j\ell}\left( \K - \frac{\p}{2}   \right) K^\ast = 0,
\end{multline}
and using that $\left(\B_+(\K),  \B_-(\K) \right)$ forms a basis of $\mathbb{R}^2$ for all $\K$, $\widehat{\W_1}^{j\ell}$ can be decomposed as 
\begin{equation}
\label{def:hatW1_cross}
\widehat{\W}_1^{j\ell} (\p,\K)=\sum_{i,m=\pm}\alpha_{im}^{j\ell} (\x,\p)\B_i\left(\K+ \frac{\p}{2}\right) \B_m^*\left(\K- \frac{\p}{2}\right).
\end{equation}
Projecting Eq.~\eqref{eq:fourier_first_order_corrector} on the left on $\C_i^*\left(\K+\frac{\p}{2}\right)$, and on the right on $\C_m\left(\K-\frac{\p}{2}\right)$, we obtain
\begin{multline}
\label{def:alpha_cross}
\alpha_{im}^{j\ell} (\x,\p)=\\
\frac{1}{2c_0^2} \left(\frac{\widehat{\nu}_1(\x,\p) \lambda_i\left( \K + (\p/2)   \right)a_m^{j\ell} \left( \K - (\p/2)   \right)   -   \widehat{\nu}_2(\x,\p)      \lambda_m\left( \K - (\p/2)   \right)  a_i^{j\ell}\left( \K + (\p/2)   \right)     }{ \lambda_i \left( \K + (\p/2)   \right) - \lambda_m\left( \K - (\p/2)   \right)    +\theta  }\right).
\end{multline} 
where we have used the three following relations 
\begin{equation*}
\lambda_\pm^* = - \lambda_\pm, \qquad \B_m^*(\p)K^*\C_i(\q)=\frac{1}{2c_0^2}\lambda_i(\q),\qquad \text{and}\qquad \C_i^*(\p)K\B_m(\q)=-\frac{1}{2c_0^2}\lambda_i(\p).
\end{equation*}

\subsection{Derivation of the transport equation}
\label{sec:derivation_RTE}

To conclude and derive the transport equation, we have to discuss the second order term $\W_2^{j\ell}$ of Eq.~\eqref{eq:asymptotic_expansion_1_annex_B}. Equating the terms associated to powers of $\epsilon$ that appears in Eq.~\eqref{eq:wigner_equation_rescale_annex_B}, through Eq. \eqref{eq:asymptotic_expansion_1_annex_B}  and Eq.~\eqref{eq:asymptotic_expansion_0_annex}, we obtain
\begin{multline}
\label{eq:second_order_cross_annex}
\Pdiff_0\left( \I \K + \frac{\mathbf{D}_\y}{2}   \right)\W_2^{j\ell} +\W_2^{j\ell} \Pdiff_0^*\left( \I \K - \frac{\mathbf{D}_\y}{2}   \right) + \mathcal{K}^{1-} K \W_1^{j\ell}+  \mathcal{K}^{2+} \W_1^{j\ell} K^* \\
+ \partial_t \W_0^{j\ell} + \Pdiff_1(\I\K)\W_0^{j\ell} + \W_0^{j\ell} \Pdiff_1^*(\I \K) =0,
\end{multline}
where the last terms do not depend on $\mathbf{D}_\y$ because $\W_0^{j\ell}$ does not depend on the $\y$-variable. Thanks to the decomposition~\eqref{eq:dec_W0}, the term $\W_2^{j\ell}$ can be chosen as being orthogonal to $\W_0^{j\ell}$ ($\W_0^{j\ell}$ is expanded over a two dimensional basis in a four dimensional vector space for any fixed $\K$) so that we necessarily have 
\[
\C^*_+(\K)\W_2^{j\ell}(t,\x,\y,\K)\C_+(\K) = 0.
\]
As a result, projecting Eq.~\eqref{eq:second_order_cross_annex} on the left on $\C^*_+(\K)$, and on the right on $\C_+(\K)$, we obtain 
\begin{equation}
\label{eq:ed0}
\partial_t a^{j\ell} + \mathcal{L}_1 \W_1^{j\ell} + \mathcal{L}_2 \W_0^{j\ell} = 0,
\end{equation}
with
\[
\mathcal{L}_1\W_1^{j\ell}(\K) :=  \C^*_+(\K)( \mathcal{K}^{1-} K \W_1^{j\ell} + \mathcal{K}^{2+} \W_1^{j\ell} K^*) \C_+(\K),
\]
and 
\[
\mathcal{L}_2\W_0^{j\ell}(\K) :=  \C^*_+(\K)\big(\Pdiff_1(\I\K)\W_0^{j\ell} + \W_0^{j\ell} \Pdiff_1^*(\I \K)\big)\C_+(\K) = c_0 \nabla_\K q(i\K)\cdot \nabla_\x a^{j\ell}= c_0 \widehat{\K} \cdot \nabla_\x a^{j\ell},
\]
remembering that $\Pdiff_1$ is defined through (\ref{eq:asymptotic_expansion_0_annex} - \ref{eq:asymptotic_expansion_annex}).

Factorizing $\widehat{\W}_1^{j\ell}$ as
\begin{equation}
\label{def:calW1}
\widehat{\W}_1^{j\ell}(\x,\p,\K)= \widehat{\nu}_1(\x,\p) \Z_1^{j\ell}(\p,\K)  +   \widehat{\nu}_2(\x,\p) \Z_2^{j\ell}(\p,\K),
\end{equation}
according to Eq.~\eqref{def:hatW1_cross} and Eq.~\eqref{def:alpha_cross}, and invoking the same mixing argument as in \cite{Messaoudi2022,Bal2005,Ryzhik1996}, yields after averaging Eq.~\eqref{eq:ed0}
\begin{equation}
\label{eq:RTE_tmp_cross}
\partial_t \langle a^{j\ell} \rangle + c_0 \widehat{\K} \cdot \nabla_\x \langle a^{j\ell} \rangle + \big\langle\mathcal{L}_1\W_1^{j\ell}\big\rangle = 0,
\end{equation}
with
\begin{align}
\label{eq:expL1}
\left< \widehat{\mathcal{L}_1\W_1^{j\ell}}(\x,\p, \K) \right>  = \int_{\mathbb{R}^3} & \left<\widehat{\nu}_1(\x,\R)\widehat{\nu}_1(\x,\p-\R)\right>K \left<\Z_1^{j\ell} \left(\p-\R,\K-\frac{\R}{2}\right) \right>\nonumber\\
& +  \left<\widehat{\nu}_1(\x,\R)\widehat{\nu}_2(\x,\p-\R)\right>K\left<\Z_2^{j\ell}\left(\p-\R,\K-\frac{\R}{2}\right)  \right>\nonumber\\
& +\left<\widehat{\nu}_2(\x,\R)\widehat{\nu}_1(\x,\p-\R)\right>\left<\Z_1^{j\ell}\left(\p-\R,\K+\frac{\R}{2}\right) \right>K^\ast\\
& + \left<\widehat{\nu}_2(\x,\R)\widehat{\nu}_2(\x,\p-\R)\right>\left<\Z_2^{j\ell}\left(\p-\R,\K+\frac{\R}{2}\right) \right>K^\ast\frac{d\R}{(2\pi)^3}.    \nonumber 
\end{align}
Moreover, according to Eq.~\eqref{def:N} and Eq.~\eqref{def:random_fluctuations_annex}, the Fourier transform $\widehat{\nu}_1$ and $\widehat{\nu}_2$ (w.r.t. the fast variable) can be written as
\begin{align}
\label{eq:fourier_transforms_fluctuations_1}
\widehat{\nu}_1(\x,\p) &= \e^{-2\I x_{0,n} p_n /\epsilon}  \e^{-\I (j+\ell)H p_n   /\epsilon} \nonumber \\
& \times\left(\sum_{m \in \mathbb{Z}} \e^{-2\I m H p_n  /\epsilon} \widehat{\nu}(\x^{j\ell}_1,\p)\mathbf{1}_ {[2mH,(2m+1)H)}(x_n- 2x_{0,n}- (j+\ell)H ) \right. \nonumber\\
& \left.  +  \sum_{m \in \mathbb{Z}} \e^{-2\I m H p_n/\epsilon} \widehat{\nu}(\x^{j\ell}_1,(\p_\perp,-p_n))\mathbf{1}_ {[2(m-1)H,2mH]}(x_n- 2x_{0,n}- (j+\ell)H)\right),
\end{align}
and
\begin{align}
\label{eq:fourier_transforms_fluctuations_2}
\widehat{\nu}_2(\x,\p) &= \e^{-\I (j+\ell)H p_n/\epsilon}  \left(\sum_{m \in \mathbb{Z}} \e^{-2\I m H p_n /\epsilon} \widehat{\nu}(\x^{j\ell}_2,\p)\mathbf{1}_ {[2mH,(2m+1)H)}(x_n- (j+\ell)H ) \right. \nonumber\\
&\left. \hspace{0.5cm} +  \sum_{m \in \mathbb{Z}} \e^{-2\I m H p_n /\epsilon} \widehat{\nu}(\x^{j\ell}_2,(\p_\perp,-p_n))\mathbf{1}_ {[2(m-1)H,2mH]}(x_n- (j+\ell)H)\right),
\end{align}
with
\[
\x^{j\ell}_1 := (\x_\perp, x_n- 2x_{0,n}- (j+\ell)H )\qquad\text{and}\qquad \x^{j\ell}_2:=(\x_\perp, x_n- (j+\ell)H).
\]
From Eq.~\eqref{eq:fourier_transforms_fluctuations_1}, Eq.~\eqref{eq:fourier_transforms_fluctuations_2}, and knowing that $\widehat{R}(\R)=\widehat{R}(\R_\perp,-\R_n)$, we obtain for the covariance function 
\begin{equation}
\label{eq:cross_power_spectra}
\left<\widehat{\nu}_1(\x,\R)\widehat{\nu}_2(\x,\p-\R) \right>=\e^{-2\I x_{0,n} r_n/\epsilon} C_1(2x_{0,n},0) + \e^{-2\I ( x_{0,n} + (j+\ell)H )r_n/\epsilon} C_2(2x_{0,n},0),
\end{equation}
where
\begin{align*}
C_1(x,x') &= (2\pi)^3  c_0^4 \delta(\p) \widehat{R}(\R)  \\
&\times \sum_{m_1,m_2\in \mathbb{Z}}  \e^{-2\I (m_1-m_2) H r_n/\epsilon} \Big( \mathbf{1}_ {[2m_1H,(2m_1+1)H)}(x_n- x- (j+\ell)H ) \\
 & \hspace{5cm}\times  \mathbf{1}_ {[2m_2H,(2m_2+1)H)}(x_n - x' -(j+\ell)H )    \nonumber  \\
& \hspace{4.5cm}+ \mathbf{1}_ {[2(m_1-1)H,2m_1H]}(x_n- x- (j+\ell)H)   \\
& \hspace{5cm} \times \mathbf{1}_ {[2(m_2-1)H,2m_2 H]}(x_n-x'-(j+\ell)H )   \Big), \nonumber
\end{align*}
and
\begin{align*}
C_2(x,x')  &= (2\pi)^3  c_0^4   \delta(\p_\perp) \delta(2r_n-p_n) \widehat{R}(\R) \\
&\times \sum_{m_1,m_2\in \mathbb{Z}}\e^{ -2\I (m_1 + m_2)H r_n /\epsilon}\Big(  \mathbf{1}_ {[2m_1H,(2m_1+1)H)}(x_n-x- (j+\ell)H ) \\
& \hspace{5cm} \times  \mathbf{1}_ {[2(m_2-1)H,2m_2H]}(x_n - x' -(j+\ell)H  )  \nonumber \\
& \hspace{4.5cm} +    \mathbf{1}_ {[2(m_1-1)H,2m_1H]}(x_n-x- (j+\ell)H  ) \\ &\hspace{5cm} \times    \mathbf{1}_ {[2m_2H,(2m_2+1)H)}(x_n-x' -(j+\ell)H  )                 \Big).\nonumber
\end{align*}
Due to the presence of the oscillatory terms $\e^{-2\I x_{0,n} r_n/\epsilon} $ and $\e^{-2 \I ( x_{0,n} + (j+\ell)H )r_n /\epsilon}$ in Eq.~\eqref{eq:cross_power_spectra}, we deduce from the Riemann-Lebesgues theorem that
\begin{equation*}
\lim\limits_{\epsilon \to 0} \int_{\mathbb{R}^3} \left<\widehat{\nu}_1(\x,\R)\widehat{\nu}_2(\x,\p-\R)\right>K \left< \Z_2^{j\ell}\left(\p-\R,\K-\frac{\R}{2}\right)\right>\frac{d\R}{(2\pi)^3}  =0,
\end{equation*}
and
\begin{equation*}
\lim\limits_{\epsilon \to 0} \int_{\mathbb{R}^3} \left<\widehat{\nu}_2(\x,\R)\widehat{\nu}_1(\x,\p-\R)\right>\left<\Z_1^{j\ell}\left(\p-\R,\K+\frac{\R}{2}\right)\right>K^* \frac{d\R}{(2\pi)^3}  =0.
\end{equation*}
Moreover, as for Eq.~\eqref{eq:cross_power_spectra}, we have
\[
\left<\widehat{\nu}_1(\x,\R)\widehat{\nu}_1(\x,\p-\R) \right> = C_1(2x_{0,n}, 2x_{0,n}) + C_2(2x_{0,n}, 2x_{0,n}),
\]
but without rapid phases in front of $C_1$ and $C_2$. However, the ones inside the sums in $C_1$ and $C_2$ implies that the non-negligible terms are obtained for $m_2=m_1$, and $m_2=-m_1$ respectively. However, knowing that 
\begin{equation*}
[2(m-1)H,2mH] \cap [2mH,2(m+1)H) = \emptyset \qquad  \forall m \in\mathbb{Z},
\end{equation*}
we necessarily have $C_2=0$. As a result, the same strategy holds for $\widehat{\nu}_2$, and we obtain at the limit
\begin{equation}
\label{eq:self_power_spectra}
\left<\widehat{\nu}_1(\x,\R)\widehat{\nu}_1(\x,\p-\R) \right> = \left<\widehat{\nu}_2(\x,\R)\widehat{\nu}_2(\x,\p-\R) \right> \underset{\epsilon\to 0}{\longrightarrow} (2\pi)^3c_0^4 \delta(\p) \widehat{R}(\R), \qquad \x \in \mathbb{R}^3.
\end{equation}
Now, going back to Eq.~\eqref{eq:expL1} and considering the change of variable $\R \to \K-\q$ for the term involving $K\langle \Z_1^{j\ell}\rangle$, and $\R \to \q-\K$ for the one involving $\langle \Z_2^{j\ell}\rangle K^\ast$ yields 
\begin{multline*}
\big\langle \widehat{\mathcal{L}_1\W_1^{j\ell}}(\x, \p, \K)  \big\rangle  \underset{\epsilon \to 0}{\longrightarrow}   c_0^4\delta(\p) \int_{\mathbb{R}^3}\widehat{R}(\K-\q)\left(K \left<\Z_1^{j\ell}\left(\q-\K,\frac{\K+\q}{2}\right)\right> \right.\\
 +\left. \left<\Z_2^{j\ell}\left(\K-\q,\frac{\K+\q}{2}\right)\right>  K^*\right) d\q. 
\end{multline*}
Finally, going back to the definitions of $\Z_1^{j\ell}$ and $\Z_2^{j\ell}$, given by Eq.~\eqref{def:calW1} through Eq.~\eqref{def:hatW1_cross} and~\eqref{def:alpha_cross}, and sending the regularization term $\theta$ to $0$, knowing that in the sense of distributions 
\begin{equation*}
\frac{1}{\I x + \theta} \, \underset{\theta \searrow 0}{\longrightarrow} \, \frac{1}{\I x} + \pi\delta(x),
\end{equation*}
we obtain after some lengthy algebra
\begin{equation*}
\int_{\mathbb{R}^3} \e^{\I \y \cdot \p} \big\langle \widehat{\mathcal{L}_1\W_1^{j\ell}}(\x, \p, \K)  \big\rangle \frac{d\p}{(2\pi)^3} \underset{\theta \searrow 0}{\longrightarrow} \Sigma(\K) a ^{j\ell}(\K), 
\end{equation*}
with
\begin{equation*}
\Sigma(\K)=\frac{\pi c_0^2\vert \K \vert^2}{2(2\pi)^3}\int_{\mathbb{R}^3} \widehat{R}(\K-\q) \delta(c_0(\vert \q \vert -\vert\K\vert ))d\q. 
\end{equation*}
Consequently, remembering Eq.~\eqref{eq:RTE_tmp_cross}, we finally obtain
\begin{equation*}
\partial_t \langle a^{j\ell}\rangle(t,\x,\K)  + c_0 \widehat{\K} \cdot \nabla_\x \langle a^{j\ell} \rangle (t,\x,\K)   = -\Sigma(\K) \langle a^{j\ell}\rangle(t,\x,\K),    
\end{equation*}
which corresponds to Eq.~\eqref{eq:RTE_coherent}.

\bibliographystyle{amsplain}
\bibliography{Biblio_boundary.bib}

\end{document}